%% file: NTM2019.tex
\documentclass[12pt,a4paper]{article}
\usepackage[T1]{fontenc}
\usepackage{mathptmx}
\usepackage[dvips]{graphicx}
\usepackage{psfrag}
\usepackage{amsmath,amssymb}        
\usepackage{latexsym}
\usepackage{url}

\newtheorem{theorem}{Teorema}
\def\RR{\hbox{I\kern-.2em\hbox{R}}}
\def\ds{\displaystyle}

\def\tw{.9\textwidth}

\newcommand{\eqnsection}{
   \renewcommand{\theequation}{{\thesection.\arabic{equation}}}
   \makeatletter
   \csname @addtoreset\endcsname{equation}{section}
   \makeatother}
\eqnsection
\title{The Iterative Transformation Method}
\author{Riccardo Fazio  \\
Department of Mathematics, Computer Science,\\ Physical Sciences and Earth Sciences,\\
University of Messina,\\
Viale F. Stagno D'Alcontres 31, 98166 Messina, Italy. \\ e-mail: rfazio@unime.it \\ home-page: http://mat521.unime.it/fazio.}
\date{\today}
\begin{document}               
\maketitle

\begin{abstract}
In a transformation method, the numerical solution of a given boundary value problem is obtained by solving one or more related initial value problems.
Therefore, a transformation method, like a shooting method, is an initial value method.
The main difference between a transformation and a shooting method is that the former is conceived and derive its formulation from the scaling invariance theory.
This paper is concerned with the application of the iterative transformation method to several problems in the boundary layer theory. 
The iterative method is an extension of the T{\"o}pfer's non-iterative algorithm developed as a simple way to solve the celebrated Blasius problem.
This iterative method provides a simple numerical test for the existence and uniqueness of solutions.
Here we show how the method can be applied to problems with a homogeneous boundary conditions at infinity and in particular we solve the Sakiadis problem of boundary layer theory. 
Moreover, we show how to couple our method with Newton's root-finder.
The obtained numerical results compare well with those available in the literature.
The main aim here is that any method developed for the Blasius, or the Sakiadis, problem might be extended to more challenging or interesting problems.
In this context, the iterative transformation method has been recently applied to compute the normal and reverse flow solutions of Stewartson for the Falkner-Skan model [Comput. \& Fluids, {\bf 73} (2013) pp. 202-209].
\end{abstract}



\input{Blasius}

\section{Main aim}
Here we show how the original treatment of the Blasius problem due to T{\"o}pfer can be extended to more complex problems of boundary layer theory.
In particular, our main concern is to solve numerically the Blasius problem, and similar problems in boundary layer theory, by initial value methods derived within scaling invariance theory.
As pointed out by NA \cite[Chapters 7-9]{Na:1979:CME}, usually a given, even simple, extension of the Blasius problem cannot be solved by T\"opfer algorithm.
Therefore, in order to extend the applicability of this non-ITM an iterative
version has been developed in \cite{Fazio:1994:FSE,Fazio:1996:NAN,Fazio:1997:NTE,Fazio:1998:SAN}.
Finally, the iterative extension of T{\"o}pfer algorithm has been applied to several problems of interest: 
free boundary problems \cite{Fazio:1990:SNA,Fazio:1997:NTE,Fazio:1998:SAN},
a moving boundary hyperbolic problem \cite{Fazio:1992:MBH},
Homann and Hiemenz problems governed by the Falkner-Skan equation in \cite{Fazio:1994:FSE},
one-dimensional parabolic moving boundary problems \cite{Fazio:2001:ITM}, two variants of the Blasius problem \cite{Fazio:2009:NTM}, namely: a boundary layer problem over moving surfaces, studied first by Klemp and Acrivos \cite{Klemp:1972:MBL}, and a boundary layer problem with slip boundary condition, that has found application in the study of gas and liquid flows at the micro-scale regime \cite{Gad-el-Hak:1999:FMM,Martin:2001:BBL}, parabolic problems on unbounded domains \cite{Fazio:2010:MBF} and,  recently, see \cite{Fazio:2015:ITM}, a further variant of the Blasius problem in boundary layer theory: the so-called Sakiadis problem \cite{Sakiadis:1961:BLBa,Sakiadis:1961:BLBb}.

\section{The ITM}
The applicability of a non-ITM to the Blasius problem is a consequence of its invariance with respect to the transformation (\ref{eq:scalinv:Blasius}); note that the asymptotic boundary condition is not invariant.
Several problems in boundary-layer theory lack this kind of invariance and cannot be solved by non-ITMs
\cite[Chapters 7-9]{Na:1979:CME}. 
To overcome this drawback, we can modify the problem at hand by introducing a numerical parameter $ h $,
and require the invariance of the modified problem with respect to an extended scaling group involving $ h $, see \cite{Fazio:1996:NAN,Fazio:1997:NTE} for details.

Let us consider the class of BVPs defined by
\begin{align}\label{eq:class}
&{\displaystyle \frac{d^3 f}{d \eta^3}} = \phi\left(\eta, f,
{\displaystyle \frac{df}{d\eta}, \frac{d^{2}f}{d\eta^2}}\right) \nonumber \\[-1.5ex]
& \\[-1.5ex]
& f(0) = a \, \qquad {\displaystyle \frac{df}{d\eta}}(0) = b + c {\displaystyle \frac{d^2f}{d\eta^2}}(0) \ , \qquad
{\displaystyle \frac{df}{d\eta}}(\eta) \rightarrow d \quad \mbox{as}
\quad \eta \rightarrow \infty \ , \nonumber
\end{align}
where $a$, $b$, $c$ and $d$ are given constants.
We turn now to define an initial value method for the class of problems (\ref{eq:class}). 
To this end, we consider an embedding parameter $h$ and the extended class of problems
\begin{align}\label{eq:class:mod}
&{\displaystyle \frac{d^3 f}{d \eta^3}} = h^{(1-3 \delta)/\sigma}\phi\left(h^{-(\delta/\sigma)} \eta, h^{-1/\sigma} f,
{\displaystyle h^{(\delta-1)/\sigma} \frac{df}{d\eta}, h^{(2 \delta-1)/\sigma} \frac{d^{2}f}{d\eta^2}}\right) \nonumber \\
&f(0) = h^{1/\sigma} a \quad  {\displaystyle \frac{df}{d\eta}}(0) = h^{(1-\delta)/\sigma} b + h^{(-\delta)/\sigma} c {\displaystyle \frac{d^2f}{d\eta^2}}(0)\ , \\
&{\displaystyle \frac{df}{d\eta}}(\eta) \rightarrow d \quad \mbox{as}
\quad \eta \rightarrow \infty \ . \nonumber
\end{align}
Let us remark here that (\ref{eq:class}) is recovered from (\ref{eq:class:mod}) by setting $h=1$.
The extended problems (\ref{eq:class:mod}) are invariant, but if $d\ne 0$ the asymptotic boundary condition is not invariant, with respect to the extended scaling group of transformations
\begin{equation}\label{eq:scaling}
f^* = \lambda f \ , \qquad \eta^* = \lambda^{\delta} \eta \ , \qquad 
h^* = \lambda^{\sigma} h \ ,   
\end{equation}
with $\delta \ne 1$ and $\sigma \ne 0$.

We have to consider now the auxiliary IVP
\begin{align}\label{eq:IVP3}
&\! \! \! \! \! \! \! \! \! \! \! \! {\displaystyle \frac{d^3 f^*}{d \eta^{*3}}} = h^{*(1-3 \delta)/\sigma}\phi\left(h^{*-(\delta/\sigma)} \eta, h^{*-1/\sigma} f, {\displaystyle h^{*(\delta-1)/\sigma} \frac{df^*}{d\eta^*}, h^{*(2 \delta-1)/\sigma} \frac{d^{2}f^*}{d\eta^{*2}}}\right) \nonumber \\[-1ex]
& \\[-1.5ex]
&\! \! \! \! \! \! \! \! \! \! \! \! f^*(0) = h^{*1/\sigma} a \ , \quad  {\displaystyle \frac{df^*}{d\eta^*}}(0) = h^{*(1-\delta)/\sigma} b + h^{*(-\delta)/\sigma} c \; p \ , \quad
{\displaystyle \frac{d^2f^*}{d\eta^{*2}}}(0) = p \ , \nonumber
\end{align}
where $p$ is defined by the user.
We have to solve (\ref{eq:IVP3}) in $[0, \eta_{\infty}^*]$, where $\eta_{\infty}^*$ is a suitable truncated boundary chosen under the asymptotic condition 
\begin{equation}
{\displaystyle \frac{df}{d\eta}}(\eta) \rightarrow d \quad \mbox{as}
\quad \eta \rightarrow \infty \ .
\end{equation} 
Now, if $d \ne 0$, then $\lambda$ is given by 
\begin{equation}\label{eq:lambda2}
\lambda = \left[ \frac{\frac{d f^*}{d \eta^{*}}(\eta^*_{\infty})}{d} \right]^{1/(1-\delta)} \ .   
\end{equation} 
Let us remark that we are able now to dismiss the above request on $d$.
In fact, if $d = 0$, then we can substitute do $d$ the value $1-h^{*(1-\delta)/\sigma}$ in (\ref{eq:class:mod}) and compute $\lambda$ by
\begin{equation}\label{eq:lambda3}
\lambda = \left[ \frac{d f^*}{d \eta^{*}}(\eta_{\infty}^*)+h^{*(1-\delta)/\sigma} \right]^{1/(1-\delta)} \ .   
\end{equation} 

It is evident that setting $h^*$ arbitrarily the transformed value of $h$ under (\ref{eq:scaling2}) can be different from one, the target value.
Therefore, we can apply a root-finder method; generally, we use the secant method.
By starting with suitable values of $ h^{*}_{0} $ and $ h^{*}_{1} $ a root-finder method is used to define the sequences $ h^{*}_{j} $ and $\lambda_j$ for $ j = 2, 3, \dots, \  $ .
A related sequence $ \Gamma (h^{*}_{j}) $, for 
$ j = 0, 1, 2, \dots, \  $ can be computed by 
\begin{equation}\label{eq:Gamma3}
\Gamma (h^{*}) = \lambda^{-\sigma} h^* - 1 \ , 
\end{equation}  
where $ \Gamma $ is defined implicitly by the solution of the IVP (\ref{eq:IVP3}).
In the following we use the notation $ \Gamma_j = \Gamma (h^{*}_{j}) $.
If the computed values of $\lambda_j$ are convergent to a value of $lambda$, that is, if $|\Gamma_j|$ goes to zero, then we can apply the inverse transformation of (\ref{eq:scaling}) to compute the numerical solution of the original boundary value problem.
A convergence criterion should be enforced, and usually I use the condition
\begin{equation}\label{eq:CC}
| \Gamma_j | \le \mbox{Tol} \ ,
\end{equation}
where $\mbox{Tol}$ is a user defined tolerance.
Once again, our goal is to find the missing initial condition ${\frac{d^2 f}{d \eta^{2}}(0)}$.

We are now ready to present the iterative method of solution in the form of an algorithm.

\noindent
{\bf The iterative algorithm.}\\
1. Input $p$, $\delta$, $ h^{*}_{0} $, $ h^{*}_{1} $, $\eta_{\infty}^*$, $d$, $\mbox{Tol}$.\\
2. $j=2, 3, \dots $; {\bf repeat} through step 5 {\bf until} condition (\ref{eq:CC}) is satisfied.\\
3. Solve  (\ref{eq:IVP3}) in the starred variables on $[0, \eta_{\infty}^*]$.\\
4. If $d \ne 0$, then compute $\lambda$ by (\ref{eq:lambda2}), else by (\ref{eq:lambda3}).\\
5. Use equation (\ref{eq:Gamma3}) to get $\Gamma_j$.\\
6. Rescale according to (\ref{eq:scaling}).\\

The above algorithm defines an ITM for the numerical solution of the class of problems (\ref{eq:class}).
In the next sections, we apply the above iterative extension of T{\"o}pfer algorithm to the Sakiadis problem, slip boundary condition, moving surface and to the Falkner-Skan model. 
But, first, let us discuss the relation between the real zero of the transformation function $\Gamma(h^*)$ and the number of solutions of the considered BVP.

\section{Existence and uniqueness}
Our main result can be stated as follows: for a given BVP the existence and uniqueness question is reduced to finding the number of real zeros of the transformation function. 
This result is proved below.

\begin{theorem}\label{Th:EU}
Let us assume that $ \delta $ and $ \sigma $ are fixed and that for every value of $ h^{*} $ the auxiliary IVP (\ref{eq:IVP3}) is well posed on $ [0, \eta^{*}_{\infty}] $. 
Then, the BVP (\ref{eq:class}) has a unique solution if and only if the transformation function has a unique real zero; nonexistence (nonuniqueness) of the solution of (\ref{eq:class}) is equivalent to nonexistence of real zeros (existence of more than one real zero) of $ \Gamma (h^*) $.
\end{theorem}

\noindent
{\bf Proof by invariance considerations.} For a proof we show that there exists a one-to-one and onto correspondence between the set of solutions of  (\ref{eq:class}) and the set of real zeros of the transformation function.
Moreover, if one of the two sets is empty the other one is empty too. 
The thesis is an evident consequence of this result.

The mentioned correspondence can be defined as follows.
For every values of  $\frac{df^*}{d\eta^*}(0)$, $\delta $ and $ \sigma $ different from zero, given a solution $ f(\eta) $ of (\ref{eq:class}), which specify a particular value of $\frac{df}{d\eta}(0)$, we can associate to it the real zero of $\Gamma$ defined by 
\[ 
h^{*} = \left[{\ds \frac{\frac{df^*}{d\eta^*}(0)}{\frac{df}{d\eta}(0)}}\right]^{\sigma/(1-2\delta)} \ . 
\]
The related value of $ \lambda = h^{*1/\sigma } $, allows us to verify 
by the substitution in (\ref{eq:Gamma3}) that we have defined a real zero of the transformation function. 
 
According to the definition of the transformation function, in general to each real zero $ h^{*} $ of $\Gamma$ there is related a solution $ f^{*}(\eta^{*}) $, defined on $ [0, \eta^{*}_{\infty}] $, of the auxiliary IVP (\ref{eq:IVP3}). 
Now, the condition for $ f^{*}(\eta^{*}), h^{*} $ to be transformed by (\ref{eq:scaling}) to $ f(\eta), 1 $ 
(where $ f(\eta) $ is defined on $ [0, \eta_{\infty}] $) is that $ \lambda = h^{*1/\sigma} $. 
Due to $ \lambda = h^{*1/\sigma} $ we have $ f^* (0) = h^{*1/\sigma} f(0) $ and $ \frac{df^*}{d\eta^*} (0) =
h^{*(1-\delta)/\sigma} \frac{df}{d\eta} (0) $, so that the relation (\ref{eq:scaling}) implies that $ f(\eta) $ verifies the boundary conditions at zero in (\ref{eq:class}). 
Hence, to each real zero of $ \Gamma $ we can associate a solution of (\ref{eq:class}).
Again $ \lambda = h^{*1/\sigma } $, so that $ f(\eta) = h^{*-1/\sigma } f^{*}(h^{*-\delta/\sigma} \eta^{*}) $.

It is easily seen that by means of the relations defined above we can fix both a right and left inverse of our correspondence.
Therefore, the correspondence is one-to-one and onto.
$ \hfill \square $

\medskip

Before proceeding further, some remarks are in order. 
First, as far as initial value problems are concerned, the theory of well-posed problems is developed in detail in several classical books, see, for instance, \cite[Chapters 2, 3 and 5]{Hartman:1982:ODE}.
In particular, the continuous dependence of the solution on parameters holds true provided suitable regularity conditions on $ \phi(\eta, f, \frac{df}{d\eta}, \frac{df^2}{d^2\eta}) $ are fulfilled.
Second, if for every value of $ h^{*} $ we assume $ \lambda (h^{*}) > 0 $, then for $ \delta \neq 0 $ and each fixed value of  $ h^* $ the scaling $ [\lambda (h^{*})]^{-\delta } \eta^{*} : [0, \eta^{*}_{\infty}] \rightarrow [0, \eta_{\infty}] $ is one-to-one and onto whereas the function of $ h^* $ defined by $ [\lambda (h^{*})]^{-\sigma } h^{*} : \RR \rightarrow \RR $ may not be one-to-one for $ \sigma \neq 0 $. 
Therefore, since $ \Gamma (h^{*}) = [\lambda (h^{*})]^{-\sigma } h^{*} - 1 $, the transformation function may not be one-to-one.
Third, by studying the behaviour of the transformation function it is possible to test the existence and uniqueness question. 

\section{Blasius and Sakiadis problems}\label{Sakiadis}
Within boundary-layer theory, the model describing the steady plane flow of a fluid past a thin plate, is given by
\begin{align}\label{PDE-model2}
& {\displaystyle \frac{\partial u}{\partial x}} +
{\displaystyle \frac{\partial v}{\partial y}} = 0 \ ,  \nonumber \\[-1.2ex]
& \\[-1.2ex]
& u {\displaystyle \frac{\partial u}{\partial x}} +
v {\displaystyle \frac{\partial u}{\partial y}} = \nu
{\displaystyle \frac{\partial^2 u}{\partial y^2}} \ ,  \nonumber 
\end{align}
where the governing differential equations, namely conservation of mass and momentum, are the 
steady-state 2D Navier-Stokes equations under the boundary layer approximations: $ u \gg v $ and the flow has a very thin layer attached to the plate,
$ u $ and $ v $ are the velocity components of the fluid in the
$ x $ and $ y $ direction, and $ \nu $ is the viscosity of the fluid.
The boundary conditions for the velocity field are
\begin{align}\label{PDE_BCs-Blasius}
& u(x, 0) = v(x, 0) = 0 \ , \qquad u(0, y) = U_\infty  \ , \nonumber \\[-1.2ex]
&\\[-1.2ex]
& u(x, y) \rightarrow U_{\infty} \quad \mbox{as}
\quad y \rightarrow \infty \ , \nonumber 
\end{align}
for the Blasius flat plate flow problem \cite{Blasius:1908:GFK}, where $ U_{\infty} $ is the main-stream velocity, and 
\begin{align}\label{PDE_BCs-Sakiadis}
& u(x, 0) = U_p \ , \qquad  v(x, 0) = 0 \ , \nonumber \\[-1.2ex]
&\\[-1.2ex]
& u(x, y) \rightarrow 0 \quad \mbox{as}
\quad y \rightarrow \infty \ , \nonumber 
\end{align}
for the classical Sakiadis flat plate flow problem \cite{Sakiadis:1961:BLBa,Sakiadis:1961:BLBb}, where $ U_p $ is the plate velocity, respectively.
The boundary conditions at $ y = 0 $ are based on the assumption that neither slip nor mass
transfer are permitted at the plate whereas the remaining boundary condition means that the velocity $ v $ tends to the main-stream velocity $ U_{\infty} $ asymptotically or gives the prescribed velocity of the plate $U_p$.

Introducing a similarity variable $\eta$ and a dimensionless stream function $f(\eta)$ as
\begin{equation}\label{eq:simvar}
\eta = y \sqrt{\frac{U}{\nu x}} \ , \qquad u = U \frac{df}{d\eta} \ , \qquad v = \frac{1}{2}\sqrt{\frac{U \nu}{x}}\left(\eta \frac{df}{d\eta}-f\right) \ ,
\end{equation}
we have
\begin{equation}\label{eq:simvar1}
\frac{\partial u}{\partial x} = -\frac{U}{2}\frac{\eta}{x}\frac{d^2f}{d\eta^2} \ , \qquad \frac{\partial v}{\partial y} = \frac{U}{2}\frac{\eta}{x}\frac{d^2f}{d\eta^2}
\end{equation}
and the equation of continuity, the first equation in (\ref{PDE-model2}), is satisfied identically.

On the other hand, we get
\begin{equation}\label{eq:simvar2}
\frac{\partial u}{\partial y} = U\frac{d^2f}{d\eta^2} \sqrt{\frac{U}{\nu x}} \ , \qquad \frac{\partial^2 u}{\partial y^2} = \frac{U^2}{\nu x}\frac{d^3f}{d\eta^3} \ .
\end{equation}
Let us notice that, in the above equations $U = U_\infty$ represents Blasius flow, whereas $U = U_p$ indicates Sakiadis flow, respectively.

By inserting these expressions into the momentum equation, the second equation in (\ref{PDE-model2}), we get 
\begin{equation}\label{eq:Bla2}
{\displaystyle \frac{d^3 f}{d \eta^3}} + \frac{1}{2} f
{\displaystyle \frac{d^{2}f}{d\eta^2}} = 0 \ ,
\end{equation}
to be considered along with the transformed boundary conditions
\begin{equation}\label{eq:Blasius:BCs}
f(0) = {\displaystyle \frac{df}{d\eta}}(0) = 0 \ , \qquad
{\displaystyle \frac{df}{d\eta}}(\eta) \rightarrow 1 \quad \mbox{as}
\quad \eta \rightarrow \infty \ , \nonumber
\end{equation}
for the Blasius flow, and
\begin{equation}\label{eq:Sakiadis:BCs}
f(0) = 0 \ , \qquad {\displaystyle \frac{df}{d\eta}}(0) = 1 \ , \qquad
{\displaystyle \frac{df}{d\eta}}(\eta) \rightarrow 0 \quad \mbox{as}
\quad \eta \rightarrow \infty \ , \nonumber 
\end{equation}
for the Sakiadis flow, respectively.

Sakiadis studied the behaviour of boundary layer flow, due to a moving flat plate immersed in an otherwise quiescent fluid, \cite{Sakiadis:1961:BLBa,Sakiadis:1961:BLBb}.
He found that the wall shear is about 34\% higher for the Sakiadis flow compared to the Blasius case. 
Later, Tsou and Goldstein \cite{Tsou:1967:FHT} made an experimental and theoretical treatment of Sakiadis problem to prove that such a flow is physically realizable.

\section{Extension of T{\"o}pfer algorithm}
Within this section, we explain how it is possible to extend T{\"o}pfer algorithm to the Sakiadis problem, that we rewrite here for the reader convenience
\begin{align}\label{eq:Sakiadis:bis}
&{\displaystyle \frac{d^3 f}{d \eta^3}} + \frac{1}{2} f
{\displaystyle \frac{d^{2}f}{d\eta^2}} = 0 \nonumber \\[-1.5ex]
&\\[-1.5ex]
&f(0) = 0 \ , \qquad {\displaystyle \frac{df}{d\eta}}(0) = 1 \ , \qquad
{\displaystyle \frac{df}{d\eta}}(\eta) \rightarrow 0 \quad \mbox{as}
\quad \eta \rightarrow \infty \ . \nonumber 
\end{align}

In order to define the ITM, we introduce the extended problem
\begin{align}\label{eq:Sakiadis:mod}
&{\displaystyle \frac{d^3 f}{d \eta^3}} + \frac{1}{2} f
{\displaystyle \frac{d^{2}f}{d\eta^2}} = 0 \nonumber \\[-1.5ex]
& \\[-1.5ex]
&f(0) = 0 \ , \quad {\displaystyle \frac{df}{d\eta}}(0) = h^{1/2} \ , \quad
{\displaystyle \frac{df}{d\eta}}(\eta) \rightarrow 1 - h^{1/2} \quad \mbox{as}
\quad \eta \rightarrow \infty \ . \nonumber
\end{align}
In (\ref{eq:Sakiadis:mod}), the governing differential equation and the two initial conditions are invariant, the asymptotic boundary condition is not invariant, with respect to the extended scaling group
\begin{equation}\label{eq:scaling2}
f^* = \lambda f \ , \qquad \eta^* = \lambda^{-1} \eta \ , \qquad 
h^* = \lambda^{4} h \ .   
\end{equation}
Moreover, it is worth noticing that the extended problem (\ref{eq:Sakiadis:mod}) reduces to the Sakiadis problem  (\ref{eq:Sakiadis:bis}) for $h=1$.
So that, to find a solution to the Sakiadis problem we have to find a zero of the so-called {\it transformation function} 
\begin{equation}\label{eq:AA2.9}
\Gamma (h^{*}) = \lambda^{-4} h^* - 1 \ , 
\end{equation}  
where the group parameter $ \lambda $ is defined by the formula
\begin{equation}\label{eq:lambda}
\lambda = \left[\displaystyle \frac{df^*}{d\eta^*}(\eta_\infty^*)+{h^*}^{1/2}\right]^{1/2} \ ,
\end{equation}
and to this end we can use a root-finder method.

Let us notice that $\lambda$ and the transformation function are defined implicitly by the solution of the IVP
\begin{align}\label{eq:Sakiadis:IVP}
&{\displaystyle \frac{d^3 f^*}{d \eta^{*3}}} + \frac{1}{2} f^*
{\displaystyle \frac{d^{2}f^*}{d\eta^{*2}}} = 0 \nonumber \\[-1.5ex]
& \\[-1.5ex]
&f^*(0) = 0 \ , \quad {\displaystyle \frac{df^*}{d\eta^*}}(0) = {h^*}^{1/2}, \quad
{\displaystyle \frac{d^2f^*}{d\eta^{*2}}}(0) = \pm 1  \nonumber \ .
\end{align}

Here, several questions are of interest.
As far as the missing initial condition is concerned, are we allowed to use the value
\begin{equation}
\frac{d^2f^*}{{d\eta^*}^2} (0) = 1 \ ,
\end{equation}
suggested to T{\"o}pfer, as discussed in section \ref{Toepfer}, by a formal series solution of the Blasius problem?
Indeed, if the first derivative of $f$ is a monotone decreasing function, then the given boundary conditions in (\ref{eq:Sakiadis:bis}) indicate that the second derivative of $f(\eta)$ has to be negative and should go to zero as $\eta$ goes to infinity and this calls for a negative value of the missing initial condition.
Is the solution to the Sakiadis problem (\ref{eq:Sakiadis:bis}) unique?
By studying the behaviour of the transformation function $\Gamma$ we can answer both questions, and this is done in the next section.

\section{Numerical Results}
It is evident that our numerical method is based on the behaviour of the transformation function.
So that, our interest is to study the behaviour of this function with respect to its independent variable, as well as the involved parameters.
We notice that, because of the two terms $h^{1/2}$, which have been introduced in the modified boundary conditions in (\ref{eq:Sakiadis:mod}), we are allowed to consider only positive values of $h^*$. 

\begin{figure}[!hbt]
	\centering
\psfrag{h*}[][]{$ h^* $} 
\psfrag{G}[][]{$ \Gamma(h^*) $} 
\psfrag{df}[][]{${\displaystyle \frac{d^2f^*}{{d\eta^*}^2}} (0) = 1$} 
\includegraphics[width=\tw]{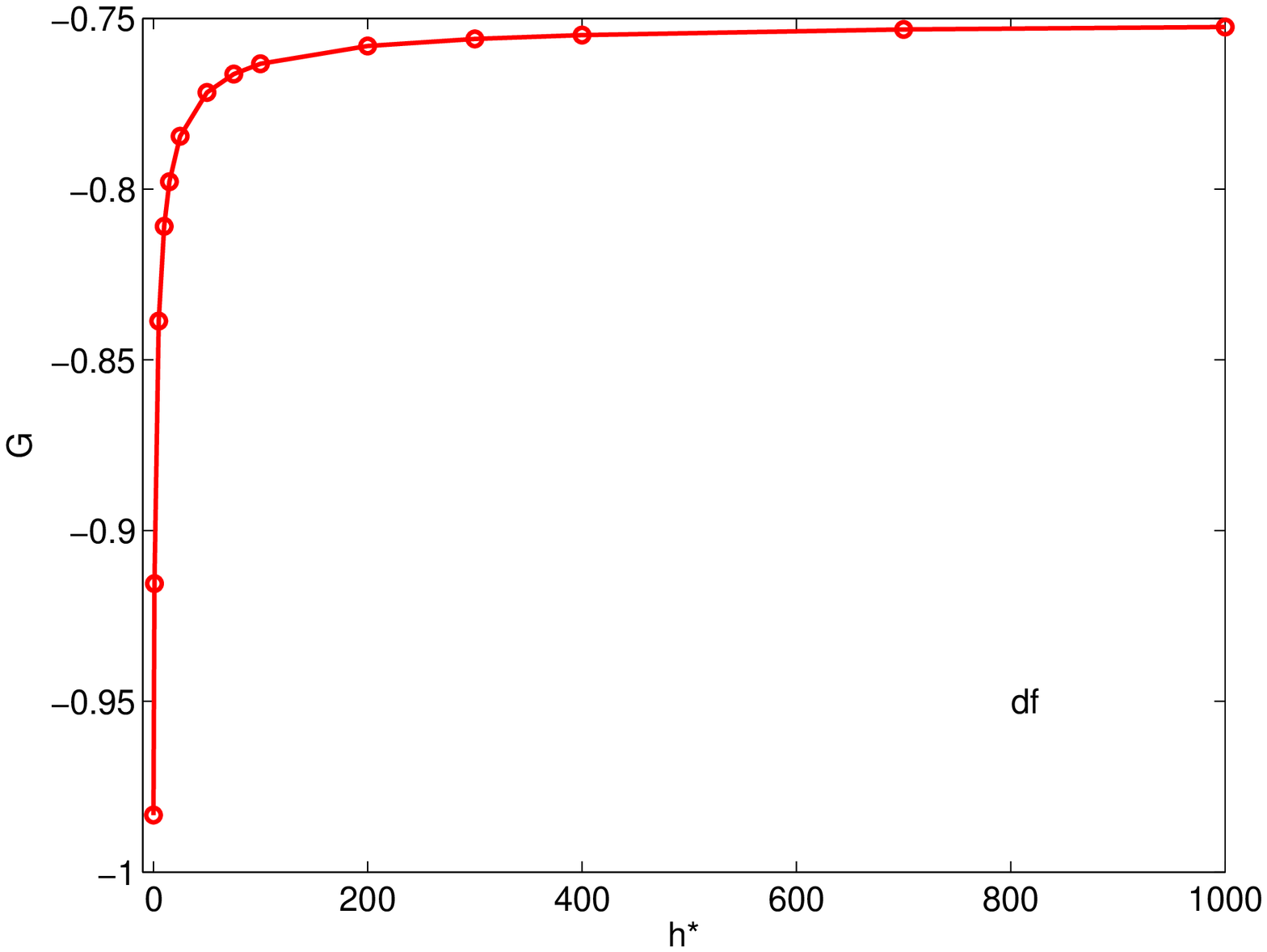} 
\caption{Plot of $\Gamma(h^*)$ for ${\displaystyle \frac{d^2f^*}{{d\eta^*}^2}}(0) = 1$.} 
	\label{fig:Gamma1}
\end{figure}

\begin{figure}[!hbt]
	\centering
\psfrag{h*}[][]{$ h^* $} 
\psfrag{G}[][]{$ \Gamma(h^*) $} 
\psfrag{df}[r][]{${\displaystyle \frac{d^2f^*}{{d\eta^*}^2}} (0) = -1$} 
\includegraphics[width=.9\textwidth]{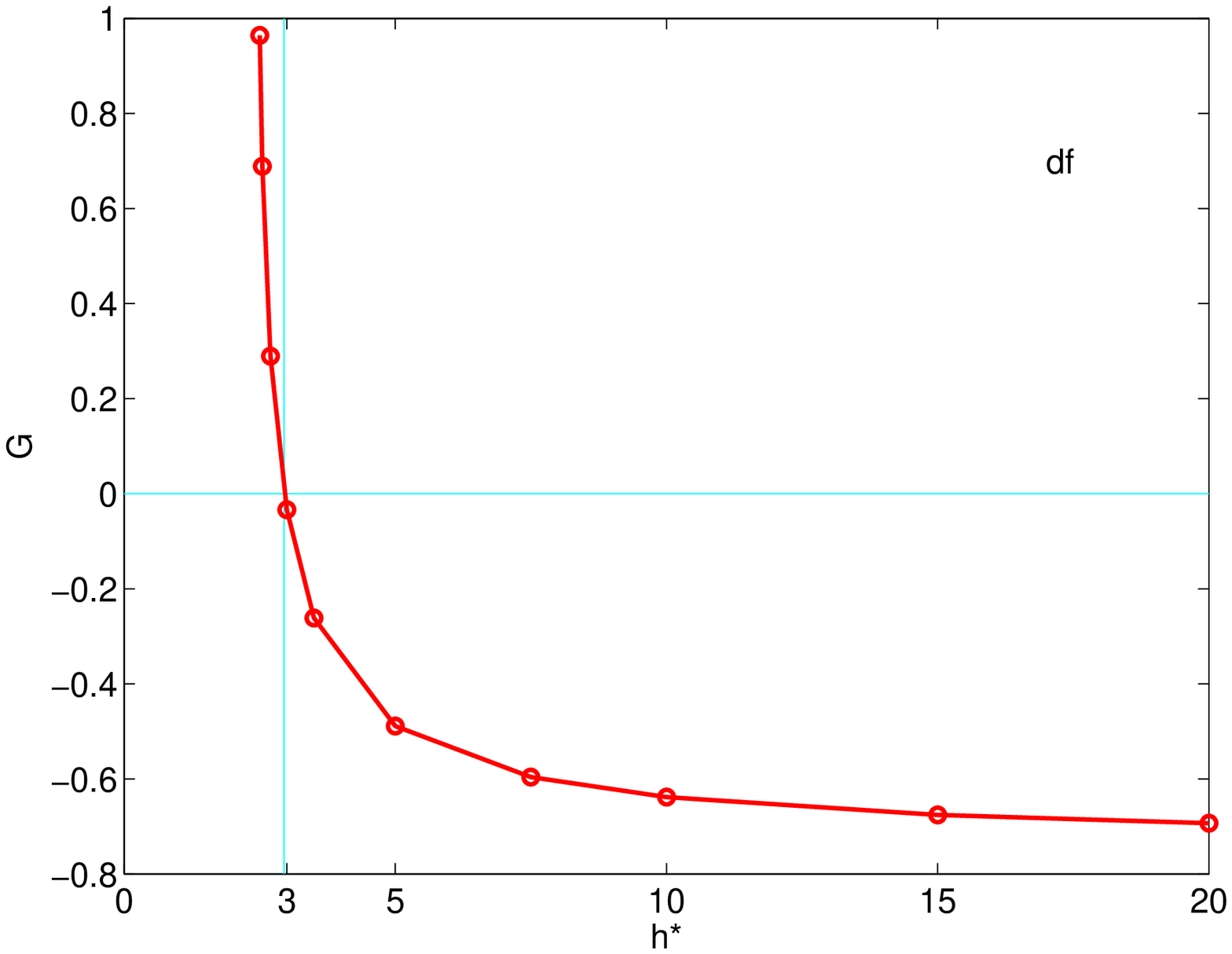} 
\caption{Plot of $\Gamma(h^*)$ for ${\displaystyle \frac{d^2f^*}{{d\eta^*}^2}} (0) = -1$.} 
	\label{fig:Gamma2}
\end{figure}

Figures \ref{fig:Gamma1}-\ref{fig:Gamma2} show the results of our numerical study concerning the dependence of $\Gamma$ with respect to the missing initial condition $\frac{d^2f^*}{{d\eta^*}^2} (0)$.
From figure \ref{fig:Gamma1} we realize that the missing initial condition cannot be positive.

For a negative missing initial condition, the numerical results are shown on figure \ref{fig:Gamma2}.
It is evident from figure \ref{fig:Gamma2} that the transformation function has only one zero and, by the theorem \ref{Th:EU} -- see also \cite{Fazio:1997:NTE} --, this means that the considered problem has one and only one solution.
Moreover, we remark that the tangent to the $\Gamma$ function at its unique zero and the $h^*$ axis defines a large angle.
This means that the quest for the $h^*$ corresponding to $h=1$ is a well-conditioned problem.
 
For a problem, in boundary layer theory, admitting more that one solution or none, depending on the value of a parameter involved see the next sections or \cite{Fazio:2009:NTM,Fazio:2013:BPF}.

For the numerical results reported in this section, the ITM was applied by setting the truncated boundary $ \eta_\infty^* = 10$.
Moreover, these results were obtained by an adaptive fourth-order Runge-Kutta IVP solver. 
The adaptive solver uses a relative and an absolute error tolerance, for each component of the numerical solution, both equal to $1 \mbox{D}{-06}$. 
Here and in the following the notation $\mbox{D}-k = 10^{-k}$ means a double precision arithmetic.

\subsection{Secant root-finder}
As a first case, the initial value solver was coupled with the simple secant root-finder with a convergence criterion given by
\begin{equation}\label{eq:conv}
\left|\Gamma(h^*)\right| \le 1 \mbox{D}{-09} \ .
\end{equation}

\begin{figure}[!hbt]
	\centering
\psfrag{e}[][]{$ \eta $} 
\psfrag{f}[][]{$f$} 
\psfrag{df}[][]{$ {\displaystyle \frac{df}{d\eta}} $} 
\psfrag{ddf}[][]{$ {\displaystyle \frac{d^2f}{d\eta^2}} $} 
\includegraphics[width=\tw]{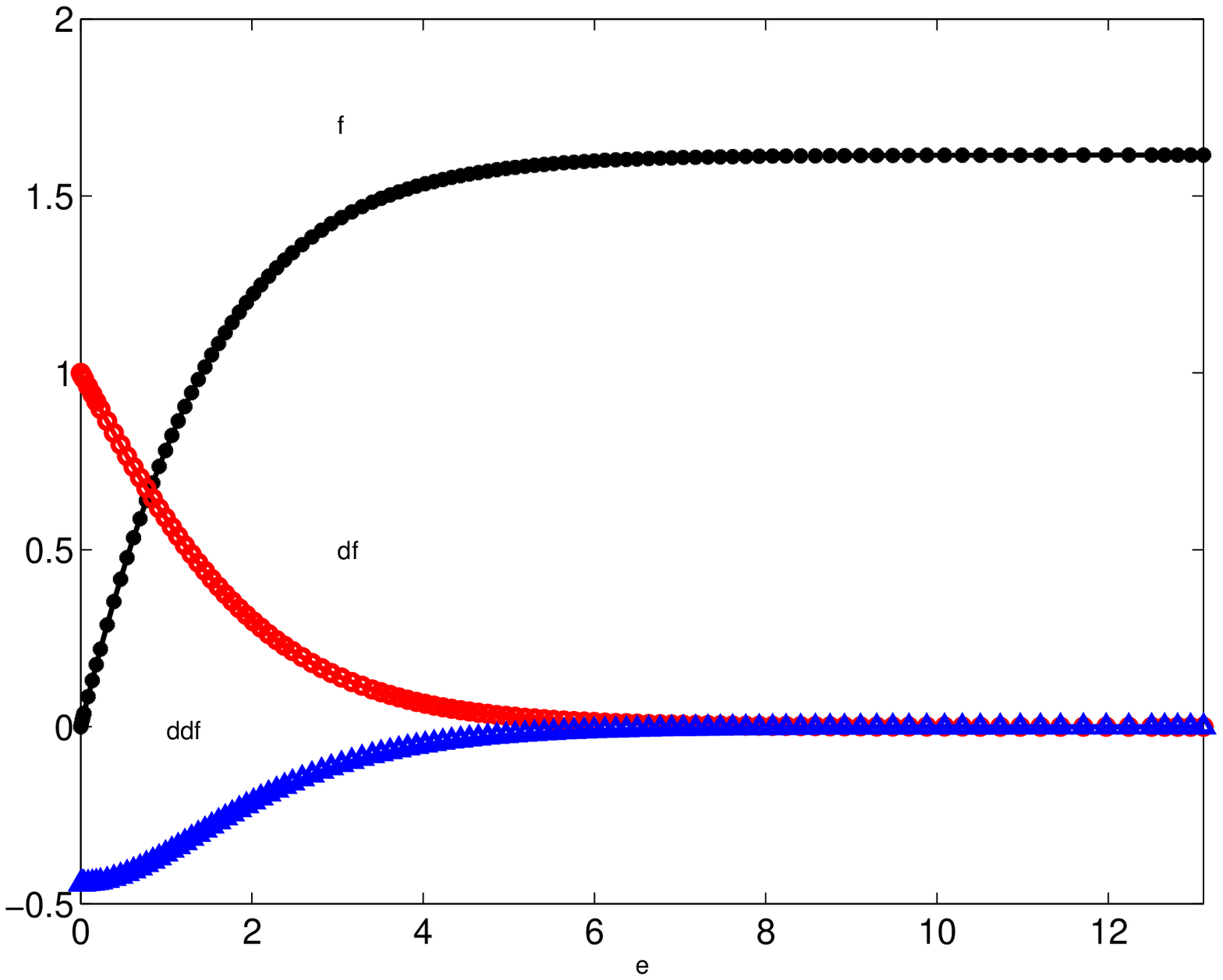} 
\caption{Sakiadis solution via the ITM.} 
	\label{fig:Sakiadis}
\end{figure}

The implementation of the secant method is straightforward.
The only difficulty we have to face is related to the choice of the initial iterates.
In this context the study of the transformation function of figure \ref{fig:Gamma2} can be helpful.
Table \ref{tab:itera} reports the iterations of our ITM.
\begin{table}[!htb]
\caption{Iterations and numerical results: secant root-finder.}
\begin{center}{ \renewcommand\arraystretch{1.3}
\begin{tabular}{rr@{.}lcr@{.}lr@{.}l} 
\hline%
{$j$}
&
\multicolumn{2}{l}%
{$h_j^*$} 
& {$\lambda_j$}
& \multicolumn{2}{c}%
{$\Gamma(h_j^*)$} 
& \multicolumn{2}{c}%
{${\displaystyle \frac{d^2f}{{d\eta}^2}} (0)$} \\[1.5ex]
\hline%
$0$  & $2$ & $5$       & $1.061732$ & $0$  & $967343$ & $-0$ & $835517$ \\
$1$  & $3$ & $5$       & $1.475487$ & $-0$ & $261541$ & $-0$ & $311310$ \\
$2$ & $3$ & $287172$ & $1.417981$ & $-0$ & $186906$ & $-0$& $350743$\\
$3$ & $2$ & $754191$ & $1.229206$ & $0$ & $206411$ & $-0$ & $538426$ \\
$4$ & $3$ & $033897$ & $1.339089$ & $-0$ & $056455$ & $-0$ & $416458$ \\
$5$ & $2$ & $973826$ & $1.318081$ & $-0$ & $014749$ & $-0$ & $436690$ \\
$6$ & $2$ & $952581$ & $1.310382$ & $0$ & $001407$ & $-0$ & $444433$ \\
$7$ & $2$ & $954432$ & $1.311058$ & $-3$ & $23\mbox{D}{-05}$ & $-0$ & $443745$ \\
$8$ & $2$ & $954391$ & $1.311043$ & $-6$ & $93\mbox{D}{-08}$ & $-0$ & $443761$ \\
$9$ & $2$ & $954391$ & $1.311043$ & $3$ & $42\mbox{D}{-12}$ & $-0$ & $443761$ \\
\hline
\end{tabular}}
\label{tab:itera}
\end{center}
\end{table}
Figure \ref{fig:Sakiadis} shows the results of our numerical approximation.
This solution was computed by rescaling, so that  ${\eta^*_\infty} < {\eta}_{\infty}$, where the chosen truncated boundary was ${\eta^*_\infty} = 10$ in our case.

\subsection{Newton's root-finder}
The same ITM can be applied by using the Newton's root-finder.
This requires a more complex treatment involving a system of six differential equations. 
Let us introduce the auxiliary variables $u_j(\eta)$ for $j = 1, 2, \dots , 6$ defined by
\begin{align}\label{eq:var}
&u_1 = f \ , \quad u_2 = \frac{df}{d\eta} \ , \quad u_3 = \frac{d^2f}{d\eta^2} \ , \nonumber \\[-1ex]
& \\[-1ex]
&u_4 = \frac{\partial u_1}{\partial h} \ , \quad u_5 = \frac{\partial u_2}{\partial h} \ , \quad u_6 = \frac{\partial u_3}{\partial h} \ . \nonumber
\end{align}
Now, the related IVP is given by
\begin{align}\label{eq:IVP6}
& \frac{du_1^*}{d\eta^*} = u_2^* \ , \nonumber \\
& \frac{du_2^*}{d\eta^*} = u_3^* \ , \nonumber \\
& \frac{du_3^*}{d\eta^*} = -\frac{1}{2}u_1^* u_3^* \ , \nonumber \\
& \frac{du_4^*}{d\eta^*} = u_5^* \ , \nonumber \\[-1.5ex]
&\\[-1.5ex]
& \frac{du_5^*}{d\eta^*} = u_6^* \ , \nonumber \\
& \frac{du_6^*}{d\eta^*} = -\frac{1}{2} \left(u_4^* u_3^* + u_1^* u_6^* \right) \ , \nonumber \\
& u_1^*(0) = 0 \ , \quad u_2^*(0) = h^{*1/2} \ , \quad u_3^*(0) = -1 \ , \nonumber \\
& u_4^*(0) = 0 \ , \quad u_5^*(0) = \frac{1}{2} h^{*-1/2} \ , \quad u_6^*(0) = 0 \ . \nonumber
\end{align}

In order to apply the Newton's root-finder, at each iteration, we have to compute the derivative
with respect to $h^*$ of the transformation function $\Gamma(h^*)$.
In our case, replacing equation (\ref{eq:lambda}) into (\ref{eq:AA2.9}), the transformation function is given by 
\begin{equation}\label{eq:Gamma2}
\Gamma(h^*) = \left[u_2^*(\eta_\infty^*)+h^{*1/2}\right]^{-2} h^* -1 \ ,
\end{equation}
and its first derivative can be easily computed as
\begin{multline}\label{eq:dGamma}
\frac{d\Gamma}{dh^*}(h^*) = \left[u_2^*(\eta_\infty^*)+h^{*1/2}\right]^{-2}\bigg\{1-2 \\ 
\left[u_5^*(\eta_\infty^*)+\frac{1}{2}h^{*-1/2}\right]\left[u_2^*(\eta_\infty^*)+h^{*1/2}\right]^{-1} h^*\bigg\}  \ .
\end{multline}
The convergence criterion is again given by (\ref{eq:conv}).
Table \ref{tab:itera2} reports the iterations of our ITM.
\begin{table}[!htb]
\caption{Iterations and numerical results: Newton's root-finder.}
\begin{center}{ \renewcommand\arraystretch{1.3}
\begin{tabular}{rr@{.}lcr@{.}lr@{.}l} 
\hline%
{$j$}
&
\multicolumn{2}{l}%
{$h_j^*$} 
& {$\lambda_j$}
& \multicolumn{2}{c}%
{$\Gamma(h_j^*)$} 
& \multicolumn{2}{c}%
{${\displaystyle \frac{d^2f}{{d\eta}^2}} (0)$} \\[1.5ex]
\hline%
$0$  & $2$ &$5$           & $1.061732$ & $0$  & $967345$ & $-0$ & $835517$ \\
$1$  & $2$ & $634888$ & $1.166846$ & $0$ & $421371$ & $-0$ & $629447$ \\
$2$ & $2$ & $812401$ & $1.255130$ & $0$ & $133241$ & $-0$& $505747$\\
$3$ & $2$ & $929233$ & $1.301740$ & $0$ & $020134$ & $-0$ & $453344$ \\
$4$ & $2$ & $953635$ & $1.310767$ & $5$ & $88\mbox{D}{-04}$ & $-0$ & $444042$ \\
$5$ & $2$ & $954391$ & $1.311043$ & $5$ & $26\mbox{D}{-07}$ & $-0$ & $443761$ \\
$6$ & $2$ & $954391$ & $1.311043$ & $4$ & $36\mbox{D}{-13}$ & $-0$ & $443761$ \\
\hline
\end{tabular}}
\label{tab:itera2}
\end{center}
\end{table}

As can be easily seen we get the same numerical results already obtained by the secant method but with a smaller numbers of iterations, cf. iteration 5 in table \ref{tab:itera2} with iteration 8 in table \ref{tab:itera}.

\section{Slip flow condition}
Let us consider again the slip boundary condition 
\begin{equation}
\frac{df}{d\eta} (0) = c \; \frac{d^2f}{d\eta^2} (0) \ ,
\end{equation}
of the previous section. 
If we need the solution for a specific value of $c$, then we can use the ITM.

In this section, we apply the ITM, and we consider a modified problem with the boundary condition
\begin{equation}
\frac{df}{d\eta} (0) = h \; c \; \frac{d^2f}{d\eta^2} (0) \ ,
\end{equation}
and the extended scaling group
\begin{equation}\label{eq:scalinv:slip2}
f^* = \lambda f \ , \qquad \eta^* = \lambda^{-1} \eta \ , \qquad 
h^* = \lambda^{-1} h \ .   
\end{equation}
Also in this case, $ \lambda $ is defined by equation (\ref{eq:lambda:moving}).
The iterative numerical results are reported in table \ref{tab:Piter}.
\input{Tab_Blasius_Slip_Iter}
The results listed in the last two columns of table \ref{tab:Piter} can be compared with similar results, obtained via a shooting method, shown in figure 1 of the proceedings report by Martin and Boyd \cite{Martin:2001:BBL}.

For the ITM, we always used $ h^*_0 = 0.1 $ and $ h^*_1 = 1 $ but for the case $ c = 50 $ where, in order to speed up the convergence, we set $ h^*_1 = 0.5 $. 
For the sake of brevity, we omit to report the iterations related to the results listed in table \ref{tab:Piter}. 
However, by setting again $ |\Gamma(\cdot)| < \mbox{D-06} $, as a convergence criterion, the Regula Falsi method converged within 8 iterations in all cases.
 Here and in the following the $\mbox{D}-k = 10^{-k}$ means a double precision arithmetic.

\section{Moving surfaces}
Klemp and Acrivos \cite{Klemp:1972:MBL} were the first to define 
the similarity model of a boundary layer problem over moving surfaces.
For this model the Blasius equation has to be considered along with the following boundary conditions
\begin{equation}\label{eq:MSCs}
f(0) = 0 \ , \qquad {\displaystyle \frac{df}{d\eta}}(0) = b \ , \ , \quad
{\displaystyle \frac{df}{d\eta}}(\eta) \rightarrow 1 \quad \mbox{as}
\quad \eta \rightarrow \infty \ ,
\end{equation}
where $ -b $ is the ratio of the plate velocity to the free stream velocity.
Let us remark here that, in contrast to the moving wall case considered in a previous section, in this case, the asymptotic boundary condition is the same as in the Blasius problem.
Klemp and Acrivos studied the effect of the parameter $ b $ on the boundary layer thickness.
For $ b < 0 $, two solutions exist only for $ b $ bigger than a critical value $ b_c $, as shown numerically by Hussaini and Lakin \cite{Hussaini:1986:ENS}.
These authors found a numerical value of $ b_c = -0.3541 $.
Hussaini et al. \cite{Hussaini:1987:SSB} proved the nonuniqueness and analyticity of solutions for $ b_c \le b $, and derived the lower bound $ -0.46824 $ for $ b_c $.

More recently, a modified Blasius equation, taking into account the effect of $ b $ on the boundary layer thickness, has been introduced by Allan \cite{Allan:1997:SSB}.
Moreover, Allan and Syam \cite{Allan:2005:ASN}, using an homotopy analysis method, defined an implicit relation between the wall shear stress and the moving wall parameters.
The study of these relation shows that two solutions exist when $b_c \approx 0.354\dots \le b$, one solution exists for $ b = b_c $ and no solution exists for $ b < b_c $.   

We have used the ITM in order to investigate the existence and uniqueness question for the Blasius model on a moving surface.
For the modified problem, we defined the boundary condition
\begin{equation}
{\displaystyle \frac{df}{d\eta}}(0) = h \; b
\end{equation}
and used the extended scaling group
\begin{equation}\label{eq:scalinv:moving}
f^* = \lambda f \ , \qquad \eta^* = \lambda^{-1} \eta \ , \qquad 
h^* = \lambda^2 h \ .   
\end{equation}
so that $ \lambda $ is defined by
\begin{equation}\label{eq:lambda:moving}
\lambda = \left[ \frac{d f^*}{d \eta^{*}}(\infty) \right]^{1/2} \ .   
\end{equation} 

\begin{figure}[!hbt]
\psfrag{e}[][]{$ \eta $} 
\psfrag{}[][]{$  $} 
\psfrag{f2}[][]{$ \frac{df}{d\eta} $} 
\psfrag{f3}[][]{$ \frac{d^2f}{d\eta^2} $}
\centering 
\includegraphics[width=\tw,height=6cm]{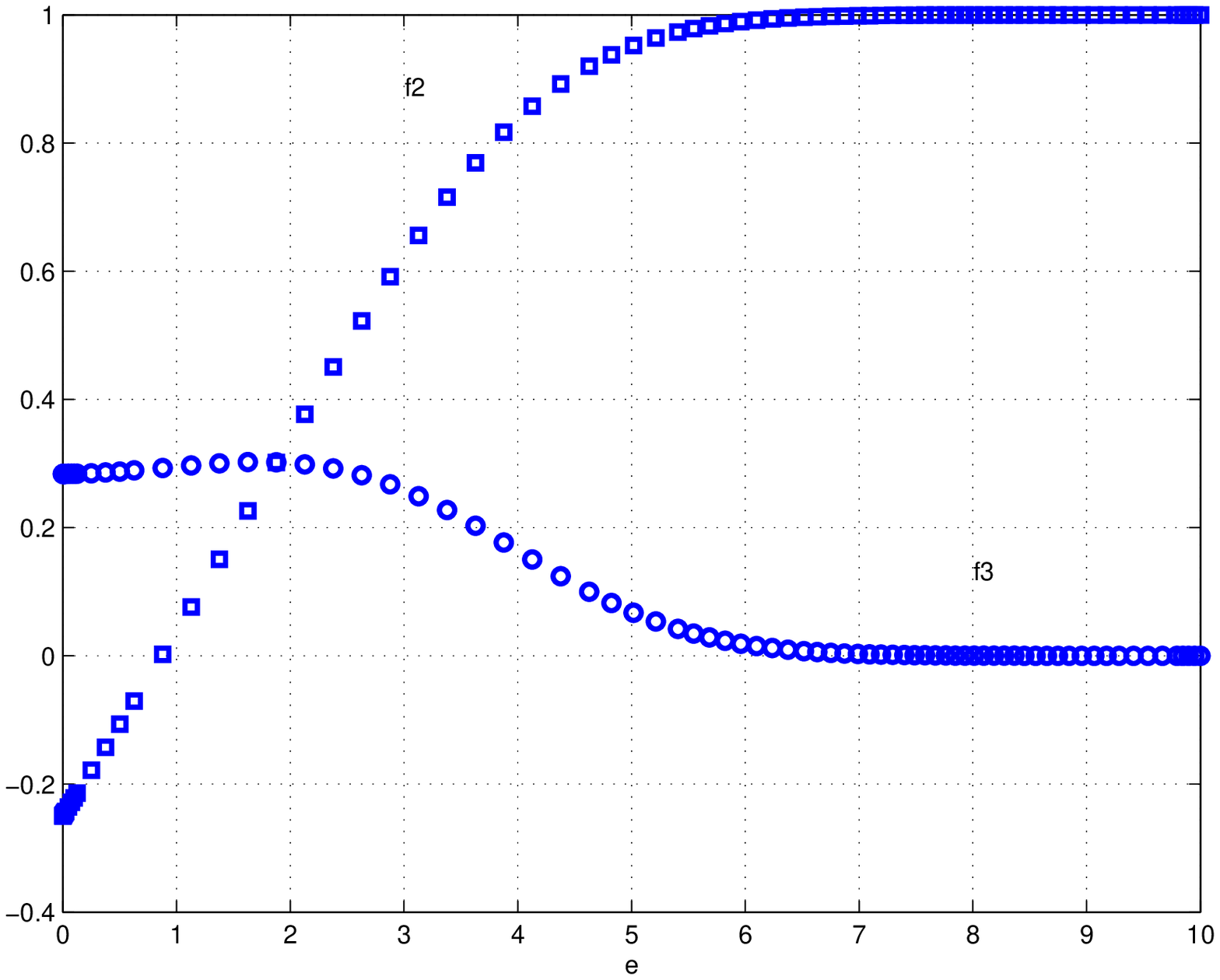} \\
\includegraphics[width=\tw,height=6cm]{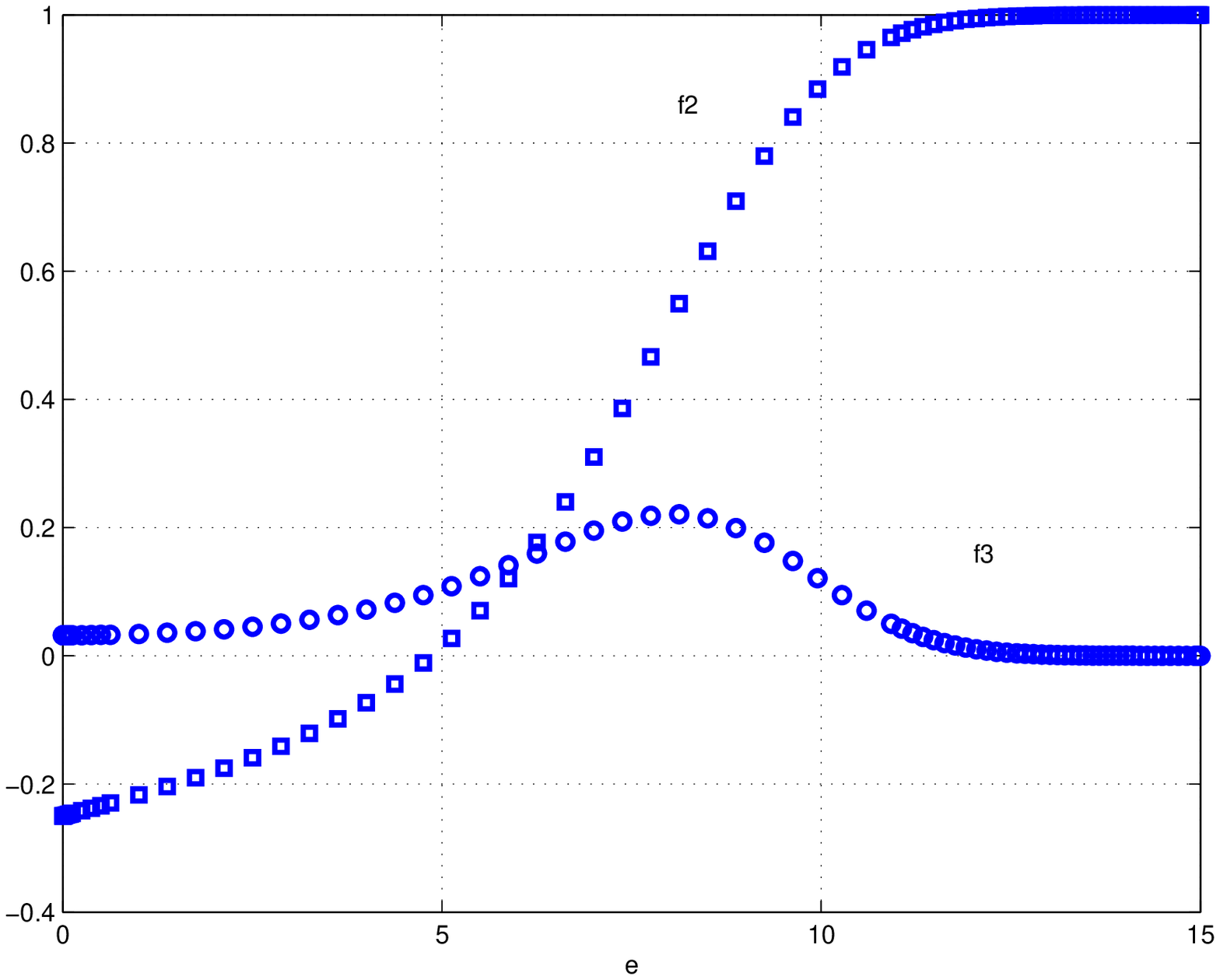}
\caption{Blasius problem with moving plate boundary conditions. The two different solutions for $ b = -0.25 $.}
	\label{fig:Blasius-Move}
\end{figure}

\input{Tab_Blasius_Move_Iter}
Let us discuss here three specific test cases.
First, we consider the case $ b = -0.25 $, and we report,  in table \ref{tab:Move}, the related numerical results found by the ITM.
In this case $ \Gamma(h^*) $ has two different zeros.
Figure \ref{fig:Blasius-Move} shows the two corresponding solutions.
It is evident, from the two frames of this figure, that the truncated boundary approach has to be supplemented by some numerical experiments, and this is more relevant in the case of nonuniqueness of the solution.
In fact, by setting $ \eta_\infty = 10 $, we miss the solution shown in the bottom frame of figure \ref{fig:Blasius-Move}. 

As a second test case, by setting $ b = -0.4 $ we find that $ \Gamma(h^*) $ has always the same negative sign, so that no solution is available for such a case.
Finally, by considering the case $ b = b_c =- 0.3541$ the ITM provided $ \frac{d^2f}{d\eta^2}(0) \approx 0.148850$.
The results obtained by the ITM in all the above cases are in agreement with the results by Hussaini and Lakin \cite{Hussaini:1986:ENS} and Allan and Syam \cite{Allan:2005:ASN}.
In particular, the behaviour of $\Gamma(h^*)$ in the second test case shows that the lower bound for $ b_c $ found by Hussaini et al. \cite{Hussaini:1987:SSB} is quite inaccurate.

For the ITM we used the Regula Falsi method as a root-finder, bracketing out the zeros of the transformation function $ \Gamma(h^*)$, along with a convergence criterion given by the inequality $ |\Gamma(h^*)| < \mbox{D-06} $.

\input{FS}

\section{Conclusions}
The applicability of a non-ITM to the Blasius problem is a consequence of the invariance of the governing differential equation and initial conditions with respect to a scaling group and the non-invariance of the asymptotic boundary condition.
Several problems in boundary-layer theory lack this kind of invariance plus non-invariance and cannot be solved by non-ITMs. 
To overcome this drawback, we can modify the problem at hand by introducing a numerical parameter $ h $, and require the invariance of the modified problem with respect to an extended scaling transformation involving $ h $, see \cite{Fazio:1996:NAN,Fazio:1997:NTE} for the application of this idea to classes of problems.
Here we show how this ITM can be used to deal with problems that admit more than one solution or with problems where the boundary condition at infinity is homogeneous. 
Moreover, we describe in details how to couple our method with Newton's root-finder.
As far as the choice of a root-finder for the ITM is concerned, we may notice that, if we limit ourselves to consider a scalar nonlinear function, then the secant method, that use one function evaluation per iteration, has an efficiency index higher than the Newton method, where we need at each iteration two function evaluations, as reported by Gautschi \cite[pp. 225-234]{Gautschi:1997:NAI}.
On the other hand, the Newton method can be preferable since it requires only one initial guess.
If we apply these methods to the solution of BVPs, jointly with a shooting or an ITM, then at each iteration the computational cost is by far higher than one or two function evaluations, and as a consequence the Newton's method might be more efficient than the secant one.


\end{document}

%% file: Blasius.tex
\section{Historical background}\label{S:History}
At the beginning of the last century Prandtl \cite{Prandtl:1904:UFK} put the foundations of
boundary-layer theory providing the basis for the unification of two, at that time seemingly incompatible, sciences: namely, theoretical hydrodynamics and hydraulics.
Boundary-layer theory has found its main application in calculating the skin-friction drag which acts on a body as it is moved through a fluid: for example the drag of an airplane wing, of a turbine blade, or a complete ship \cite{Schlichting:1979:BLT}.
Moreover, Boyd \cite{Boyd:2008:BFC} uses the problem considered by Prandtl as an example were some good analysis, before the computer invention, allowed researchers of the past to solve problems, governed by partial differential equations, that might be otherwise impossible to face.

Blasius problem \cite{Blasius:1908:GFK} is the simplest nonlinear boundary layer problem.
A study by Boyd point out how this particular problem has arisen the interest of prominent scientist, like H. Weyl, J. von Neumann, M. Van Dyke, etc., see Table 1 in \cite{Boyd:1999:BFC}. 
The main reason for this interest is due to the hope that any approach developed for this epitome can be extended to more difficult hydrodynamics problems.

Blasius main interest was to compute, without worrying about existence or uniqueness of its boundary value problem (BVP) solution -- it was Weyl who proved in \cite{Weyl:1942:DES} that Blasius problem has one and only one solution --, the value of  $\lambda$ the so-called shear stress.
To compute this value, Blasius used a formal series solution around $\eta=0$ and an asymptotic expansions for large values of $\eta$, adjusting the constant $\lambda$ so as to connect both expansions in a middle region.
In this way, Blasius obtained the (erroneous) bounds $ 0.3315 < \lambda < 0.33175$.

A few years later, T{\"o}pfer \cite{Topfer:1912:BAB} revised the work by Blasius and solved numerically the Blasius equation with suitable initial conditions and the classical order-four Runge-Kutta method.
He then arrived, without detailing his computations, at the value
$\lambda \approx 0.33206$, contradicting the bounds reported by Blasius.  

Thereafter, the quest for a good approximation of $\lambda$ was a main concern.
This is seldom the case for the most important problems of applied mathematics:
at the first study everyone would like to know if there is a method to solve a given problem, but, as soon as a problem is solved, then we would like to know how accurate is the computed solution and whether there are different methods that can provide a solution with less effort. 
By using a power series, Bairstow \cite{Bairstow:1925:SF} reports $\lambda \approx 0.335$, and Goldstein \cite{Goldstein:1930:CSS} obtains $\lambda \approx 0.332$ or, using a finite difference method, Falkner \cite{Falkner:1936:MNS} finds $\lambda \approx 0.3325765$, and Howarth \cite{Horwarth:1938:SLB} yields $\lambda \approx 0.332057$.
Fazio \cite{Fazio:1992:BPF}, using a free boundary formulation of the Blasius problem, finds $\lambda \approx 0.332057336215$.
Boyd \cite{Boyd:1999:BFC} uses T{\"o}pfer's algorithm to obtain the accurate value $\lambda \approx 0.33205733621519630$.
By the Adomain's decomposition method Abbasbandy \cite{Abbasbandy:2007:NSB} finds $\lambda \approx 0.333329$, whereas a variational iteration method with Pad\'e approximants allows Wazwaz \cite{Wazwaz:2007:VIM} to calculate, the imprecise value, $\lambda \approx 0.3732905625$.
Tajvidi et al. \cite{Tajvidi:2008:MRL} apply modified rational Legendre functions to get a value of  $\lambda \approx 0.33209$.

To compute the value of $\lambda$, we can apply also the Crocco formulation \cite{Crocco:1941:SLL}.
For instance, Vajravelu et al. \cite{Vajravelu:1991:SSS} use the Runge-Kutta method and a shooting technique to solve numerically the Crocco formulation and obtain the value $\lambda \approx 0.3322$, and
Callegari and Friedman \cite{Callegari:1968:ASN} reformulate the Blasius problem in terms of the Crocco variables, show that this problem has an analytical solution, and compute the following bounds: $0.332055 < \lambda < 0.33207$.

At the turning of this new century, as the number of applications of microelectronics devices increases, boundary-layer theory has found a renewal of interest within the study of gas and liquid flows at the micro-scale regime, see, for instance, Gad el Hak \cite{Gad-el-Hak:1999:FMM} or Martin and Boyd \cite{Martin:2001:BBL}.

Our main goal here is to show how to solve numerically Blasius problem, and similar problems in the boundary layer theory, by initial value methods derived from scaling invariance theory.
In the literature, these methods are referred to as numerical transformation methods (TM) and in recent years these methods have been effectively applied to several problems of interest.

\section{Fluid flow on a flat plate}\label{S:Fluid}
\begin{figure}[!hbt]
	\centering
\psfrag{y}[][]{$ y $} 
\psfrag{z}[][]{$ z $} 
\psfrag{O}[][]{$ O $}
\psfrag{V}[][]{$ V_\infty $} 
\includegraphics[width=\tw]{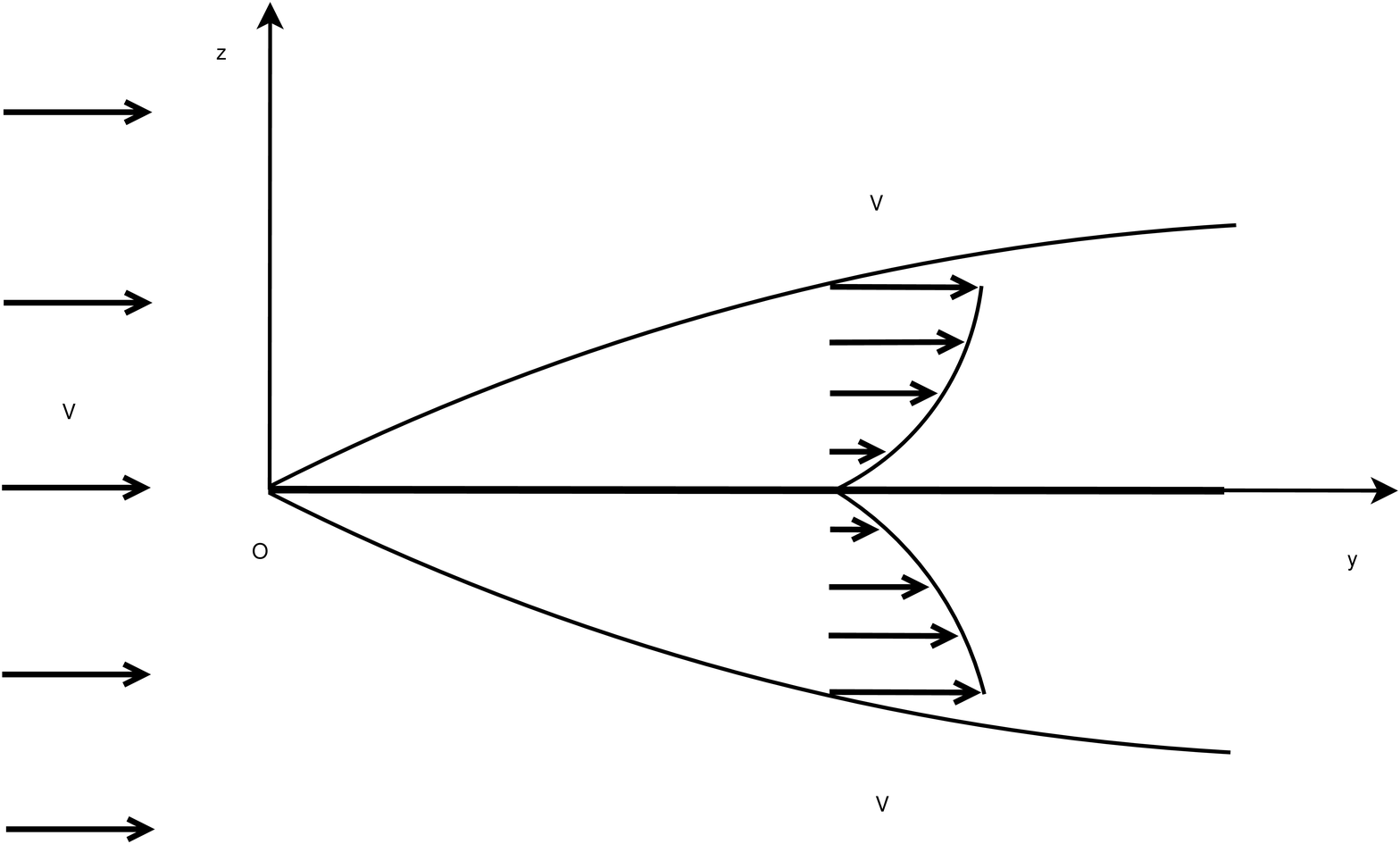}
\caption{Boundary layer over a thin plate.}
	\label{fig:plate}
\end{figure}

The model describing
the steady plane flow of a fluid past a thin plate\index{fluid!on a thin plate}%
, provided the boundary layer assumptions are verified (the flow has a very thin layer attached to the plate and $ v \gg w $), is given by
\begin{align}\label{PDE-model}
& {\displaystyle \frac{\partial v}{\partial y}} +
{\displaystyle \frac{\partial w}{\partial z}} = 0  \nonumber \\
& v {\displaystyle \frac{\partial v}{\partial y}} +
w {\displaystyle \frac{\partial v}{\partial z}} = \nu
{\displaystyle \frac{\partial^2 v}{\partial z^2}}  \\
& v(y, 0) = w(y, 0) = 0 \ , \quad v(y, z) \rightarrow V_{\infty}
\quad \mbox{as} \quad z \rightarrow \infty \ , \nonumber 
\end{align}
where the governing differential equations, namely conservation of mass and momentum, are the 
steady-state 2D Navier-Stokes
equations under the boundary layer approximations,
$ v $ and $ w $ are the velocity components of the fluid in the
$ y $ and $ z $ direction, $ V_{\infty} $
represents the main-stream velocity, see the draft in figure \ref{fig:plate}, and $ \nu $ is the viscosity of the fluid.
The boundary conditions at $ z = 0 $ are based on the
assumption that neither slip nor mass
transfer are permitted at the plate whereas the remaining boundary condition
means that the velocity $ v $ tends to the main-stream velocity
$ V_{\infty} $ asymptotically.

In order to study this problem, it is convenient to introduce a
potential (stream function) $ \psi(y, z) $ defined by
\begin{equation}
v = {\displaystyle \frac{\partial \psi}{\partial z}} \ , \qquad
w = - {\displaystyle \frac{\partial \psi}{\partial y}} \ .
\end{equation}
The physical motivation for introducing this function is that constant
$ \psi $ lines are steam-lines.
The mathematical motivation for introducing such a new variable is that the equation of continuity is satisfied identically, and we have to deal only with the transformed momentum equation.
In fact, introducing the stream function the problem can be rewritten as follows
\begin{align}\label{stream-model}
& \nu {\displaystyle \frac{\partial^3 \psi}{\partial z^3}} +
{\displaystyle \frac{\partial \psi}{\partial y}}
{\displaystyle \frac{\partial^2 \psi}{\partial z^2}} -
{\displaystyle \frac{\partial \psi}{\partial z}}
{\displaystyle \frac{\partial^2 \psi}{\partial y \partial z}}
 = 0  \nonumber \\
& {\displaystyle \frac{\partial \psi}{\partial y}}(y, 0) =
{\displaystyle \frac{\partial \psi}{\partial z}} (y, 0) = 0  \\
& {\displaystyle \frac{\partial \psi}{\partial z}} (y, z)
\rightarrow V_{\infty}  \quad \mbox{as} \quad
z \rightarrow \infty \ . \nonumber
\end{align}

\subsection{Blasius problem}
Blasius \cite{Blasius:1908:GFK} used the following similarity transformation
\begin{equation}
 \eta = z \left(\frac{V_\infty}{\nu y} \right)^{1/2} \ , \qquad  f(\eta) = \psi(y, z)
\left(\nu y V_\infty \right)^{-1/2} \ , 
\end{equation}
that reduces the partial differential model (\ref{stream-model}) to \index{Blasius problem}%
\begin{align}\label{eq:Blasius}
&{\displaystyle \frac{d^3 f}{d \eta^3}} + \frac{1}{2} f
{\displaystyle \frac{d^{2}f}{d\eta^2}} = 0 \nonumber \\[-1.5ex]
& \\[-1.5ex]
&f(0) = {\displaystyle \frac{df}{d\eta}}(0) = 0 \ , \quad
{\displaystyle \frac{df}{d\eta}}(\eta) \rightarrow 1 \quad \mbox{as}
\quad \eta \rightarrow \infty \ , \nonumber
\end{align}
i.e., a boundary value problem (BVP) defined on a semi-infinite interval.
Blasius solved this BVP by patching a power series to an asymptotic approximation at some finite value of $ \eta $.

\section{T{\"o}pfer transformation}\label{Toepfer}
In order to clarify T{\"o}pfer \cite{Topfer:1912:BAB} derivation of a further transformation of variables that reduces the BVP into an initial value problem (IVP) we consider the derivation of the series expansion solution. 
Of course, some of the coefficients of the series can be evaluated by imposing the boundary conditions at $ \eta = 0 $.
Moreover, we set
\begin{equation}
    \lambda = \frac{d^2f}{d\eta^2}(0) 
\end{equation} 
where $ \lambda $ is a nonzero constant.
So that, we look for a series solution defined as
\begin{equation}
    f(\eta) = \frac{\lambda}{2} \eta^2 + \sum_{n=3}^\infty C_n \eta^n
\end{equation}
where the coefficients $ \lambda $ and $ C_n $, for $ n = 3, 4, \dots $, are constants to be determined.
In fact, the boundary values at the plate surface, at $\eta = 0$,
require that $ C_0 = C_1 = 0 $, and we also have $ C_2 = \lambda/2 $ by the definition of $ \lambda $.
Now, we substitute this series expansion into the governing differential equation, whereupon we find 
\begin{multline} 
\sum_{n=3}^\infty n (n-1)(n-2)C_n \eta^{n-3} \\
 + \frac{1}{2} 
\left(\frac{\lambda}{2} \eta^2 + \sum_{n=3}^\infty C_n \eta^n\right)
\left[\lambda + \sum_{n=3}^\infty n (n-1) C_n \eta^{n-2} \right] = 0 
\end{multline}
or in expanded form 
\begin{multline}
 \left[3 \; 2 \; C_3 \right] 
+ \left[4 \; 3 \; 2 \; C_4 \right] \eta 
+ \left[5 \; 4 \; 3 \; C_5 + \frac{1}{2} \; 2 \; \frac{\lambda}{2} \frac{\lambda}{2}\right] \eta^2 \\
 + \left[6 \; 5 \; 4 \; C_6 + \frac{1}{2} \; 2 \; \frac{\lambda}{2} C_3 
 + \frac{1}{2} \; \frac{\lambda}{2} \; 3 \; 2 \; C_3\right] \eta^3 
 + \cdots = 0 \ .
\end{multline}
According to a standard approach, we have to require that all coefficients of the powers of $ \eta $ to be zero.
It is an easy matter to compute the coefficients of the series expansion in terms of $ \lambda $:
\begin{align*}
C_3 &= C_4 = 0 \ , \qquad C_5 = -\frac{\lambda^2}{2\; 5!} \nonumber \\
C_6 &= C_7 = 0 \ , \qquad C_8 = 11\frac{\lambda^3}{2^2\; 8!} \nonumber \\
C_9 &= C_{10} = 0 \ , \qquad C_{11} = -375\frac{\lambda^4}{2^3\; 11!} \nonumber \\
&\mbox{and so on ...} \nonumber 
\end{align*}
The solution can be written as
\begin{equation}
f = \frac{\lambda \eta^2}{2} - \frac{\lambda^2 \eta^5}{2\; 5!} 
+ \frac{11 \; \lambda^3 \eta^8}{2^2\; 8!} 
- \frac{375 \; \lambda^4 \eta^{11}}{2^3\; 11!} + \cdots 
\end{equation}
where the only unknown constant is $ \lambda $.
In principle, $ \lambda $ can be determined by imposing the boundary condition at the second point, but, in this case, this cannot be done because the left boundary condition is given at infinity.
However, by modifying the powers of $ \lambda $ we can rewrite the series expansion as
\begin{multline} 
\lambda^{-1/3} f = \frac{\left(\lambda^{1/3} \eta\right)^2}{2} - \frac{\left(\lambda^{1/3} \eta\right)^5}{2\; 5!} \\ 
+ \frac{11 \; \left(\lambda^{1/3} \eta\right)^8}{2^2\; 8!} 
- \frac{375 \; \left(\lambda^{1/3} \eta\right)^{11}}{2^3\; 11!} + \cdots 
\end{multline} 
which suggests a transformation of the form \index{scaling transformation}%
\begin{equation}\label{eq:scalinv:Blasius}
f^* = \lambda^{-1/3} f \ , \qquad \eta^* = \lambda^{1/3} \eta  \ .   
\end{equation}
In the new variables the series expansion becomes 
\begin{equation}
f^* = \frac{\eta^{*2}}{2} - \frac{\eta^{*5}}{2\; 5!} 
+ \frac{11 \; \eta^{*8}}{2^2\; 8!} 
- \frac{375 \; \eta^{*11}}{2^3\; 11!} + \cdots 
\end{equation}
which does not depend on $ \lambda $.
We notice that the governing differential equation and the initial conditions at the free surface, at $\eta = 0$, are left invariant by the new variables defined above.
Moreover, the first and second order derivatives transform in the following way
\begin{equation}
\frac{d f^*}{d \eta^{*}} = \lambda^{-2/3} \frac{d f}{d \eta} \ ,   
\qquad
\frac{d^2 f^*}{d \eta^{*2}} = \lambda^{-1} \frac{d^2 f}{d \eta^{2}} \ .  
\end{equation}
As a consequence of the definition of $ \lambda $ we have
\begin{equation}
\frac{d^2 f^*}{d \eta^{*2}} (0) = 1 \ ,   
\end{equation} 
and this explains why in these variables the series expansion does not depend on $ \lambda $.
Furthermore, the value of $ \lambda $ can be found on condition that we have an approximation for $ \frac{d f^*}{d \eta^{*}}(\infty) $, say $ \frac{d f^*}{d \eta^{*}}(\eta_{\infty}) $ where $\eta^*_{\infty}$ is a suitable truncated boundary.
In fact, by the above relation we get
\begin{equation}\label{eq:lambda:Blasius}
\lambda = \left[ \frac{d f^*}{d \eta^{*}}(\eta^*_{\infty}) \right]^{-3/2} \ .   
\end{equation} 

From a numerical viewpoint, BVPs must be solved within the computational domain simultaneously (a \lq \lq stationary\rq \rq \ problem), whereas IVPs can be solved by a stepwise procedure (an \lq \lq evolution\rq \rq \ problem).
Somehow, numerically, IVPs are easier than BVPs.  

\subsection{T{\"o}pfer algorithm}
\begin{figure}[!hbt]
	\centering
\psfrag{e}[][]{$ \eta^*, \eta $} 
\psfrag{f2star}[][]{$ \frac{df^*}{d\eta^*} $} 
\psfrag{f2}[][]{$ \frac{df}{d\eta} $} 
\psfrag{f3star}[][]{$ \frac{d^2f^*}{d\eta^{*2}} $} 
\psfrag{f3}[][]{$ \frac{d^2f}{d\eta^2} $} 
\includegraphics[width=\tw]{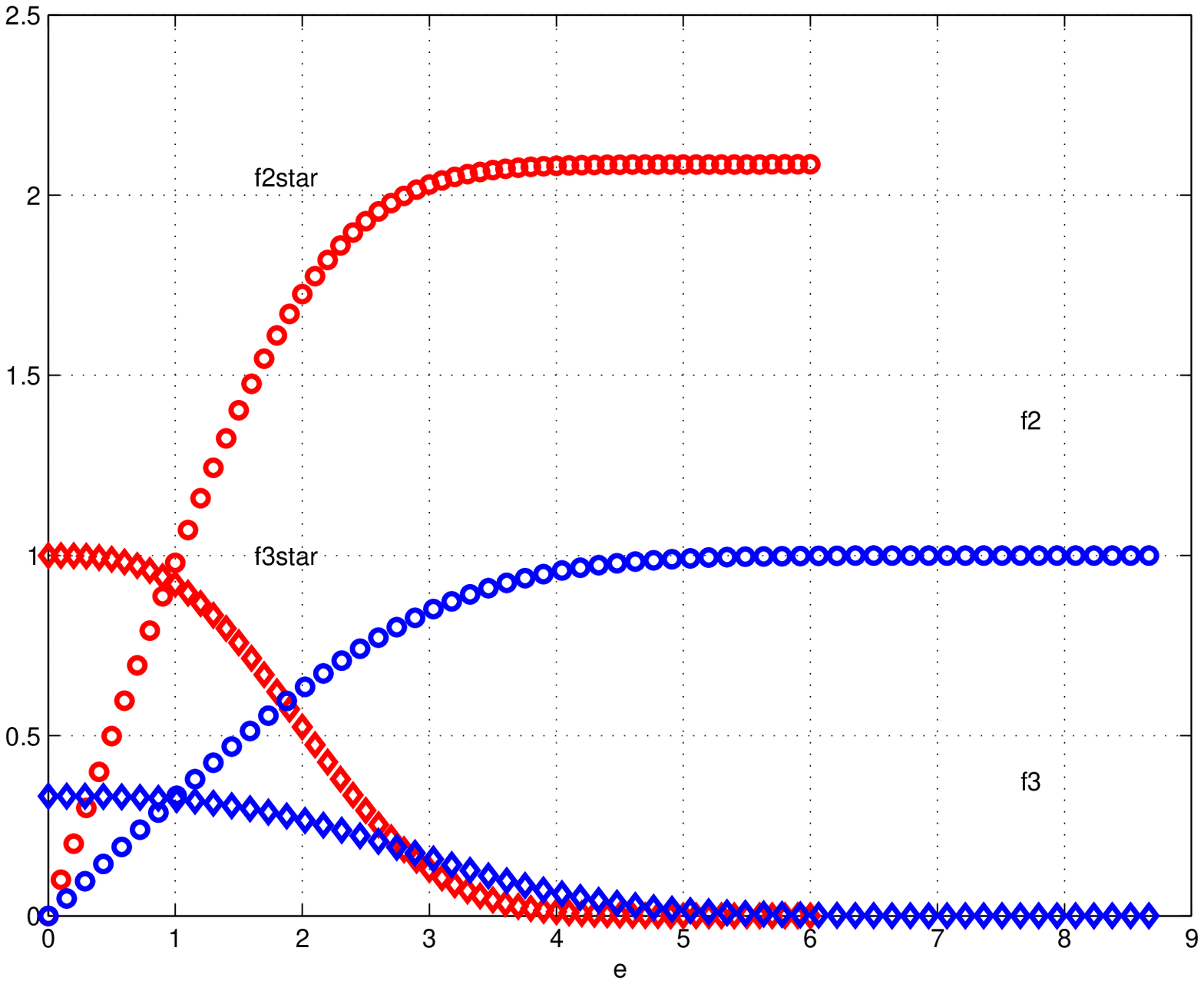}
\caption{Blasius solution by a non-ITM for $\eta^* \in [0,6]$ we get $\lambda_2 = 0.332057687$.} 
	\label{fig:Blasius2}
\end{figure}

Let us list the steps necessary to solve the Blasius problem by the T{\"o}pfer algorithm.
In this way, we define a non-iterative (I)TM.
We have to:\index{T{\"o}pfer algorithm}%
\begin{enumerate}
    \item solve the auxiliary IVP \index{IVP!auxiliary}%
\begin{align}\label{eq:Blasius2}
&{\displaystyle \frac{d^3 f*}{d \eta^{*3}}} + \frac{1}{2} f^*
{\displaystyle \frac{d^{2}f^*}{d\eta^{*2}}} = 0 \nonumber \\[-1.5ex]
& \\[-1.5ex]
&f^*(0) = {\displaystyle \frac{df^*}{d\eta^*}}(0) = 0, \qquad
{\displaystyle \frac{d^2f^*}{d\eta^{*2}}}(0) = 1  \nonumber
\end{align}
and, in particular, get an approximation for
$ \frac{d f^*}{d \eta^{*}}(\infty) $;
\item
compute $ \lambda $ by equation (\ref{eq:lambda:Blasius});
\item
obtain $ f(\eta) $, ${\displaystyle \frac{df}{d\eta}}(\eta)$, and ${\displaystyle \frac{d^2f}{d\eta^2}}(\eta)$ by the inverse transformation of (\ref{eq:scalinv:Blasius}).
\end{enumerate}
Indeed, T{\"o}pfer solved the IVP for the Blasius equation once.
At large but finite $ \eta_j^* $, ordered so that $ \eta_j^* < \eta_{j+1}^* $, we can compute by equation (\ref{eq:lambda:Blasius}) the corresponding $ \lambda_j $.  
If two subsequent values of $ \lambda_j $ agree within a specified accuracy, then $ \lambda $ is approximately equal to the common value of the $ \lambda_j $, otherwise, we can march to a larger value of $ \eta $ and try again.
Using the classical fourth order Runge-Kutta \index{Runge-Kutta method}%
 method and a grid step $ \Delta \eta^* = 0.1$
T{\"o}pfer was able to determine $ \lambda $ with an error less than $ 10^{-5} $. 
He used the two truncated boundaries $\eta_1^* = 4$ and $\eta_2^* = 6$.
In figure \ref{fig:Blasius2} we plot the more accurate numerical solution obtained by T{\"o}pfer's algorithm defined above.
We notice that this figure shows the solutions of the auxiliary IVP (\ref{eq:Blasius2}) and of the BVP (\ref{eq:Blasius}).

\section{Rubel error analysis}
The boundary condition at infinity is certainly not suitable for a numerical treatment.
This condition has usually been replaced by the same condition applied at a truncated boundary, see Collatz \cite[pp 150-151]{Collatz} or Fox \cite[p. 92]{Fox}.
In the truncated boundary formulation $ f_M(\eta) $ is defined by
\begin{eqnarray}
& {\displaystyle \frac{d^3 f_M}{d \eta^3}} + f_M
{\displaystyle \frac{d^{2}f_M}{d\eta^2}} = 0 \nonumber \\[-1.5ex]
\label{p7} \\[-1.5ex]
& f_M(0) = {\displaystyle \frac{df_M}{d\eta}}(0) = 0, \qquad
{\displaystyle \frac{df_M}{d\eta}}(M) = 1  \nonumber
\end{eqnarray}
where $ M $ represents the truncated boundary. 
It is evident that also in (\ref{p7}) the governing DE and the two boundary
conditions at the origin are left invariant by the scaling transformation (\ref{eq:scalinv:Blasius}).

For the error related to the truncated boundary solution $ f_M(\eta) $ defined by
\begin{equation}
e(\eta) = |f(\eta) - f_M(\eta)|, \qquad \eta \in [0, M] \ ,
\end{equation}
the following theorem holds true.

\begin{theorem} {\bf \ (due to Rubel \cite{Rubel:1955:EET}).} A truncated boundary formulation of the Blasius problem introduces an error which verifies the following inequality
\begin{equation}
e(\eta) \leq M {\displaystyle \frac{d^2 f_M}{d\eta^2}}(M)
[f_M(M)]^{-1}  \quad .
\end{equation}
\end{theorem}

\medskip
\noindent
{\bf Outline of the proof.} 
As proved by Weyl \cite{Weyl:1942:DES}, it is true that
\begin{equation}\label{Weylin}
{\displaystyle \frac{d^2 f}{d\eta^2}}(\eta) > 0, \qquad
\mbox{for \ \ } \eta \in (0, \infty)  \ .
\end{equation}
By (\ref{Weylin}) and taking into account the boundary conditions in (\ref{eq:Blasius}), we have that
\begin{equation}
{\displaystyle \frac{df}{d\eta}}(\eta) \ \ \mbox{and} \ \ f(\eta) \ \
\mbox{are increasing functions on} \ \ \eta \in (0, \infty) \quad .
\end{equation}
As a consequence the function $ \lambda^2 {\displaystyle \frac{df}{d\eta}}(\lambda M) $ is zero for $ \lambda = 0 $, increases with $ \lambda $ and tends to infinity as $ \lambda \rightarrow \infty $. 
For some value of $ \lambda \in (0, \infty) $ we must have
\begin{equation}
\lambda^2 {\displaystyle \frac{df}{d\eta}}(\lambda M) = 1  \quad .
\end{equation}
This value verifies $ \lambda > 1 $ because $ \lambda M $ is a
finite value, $ {\displaystyle \frac{df}{d\eta}}(\eta) $ is an
increasing function and
$ {\displaystyle \frac{df}{d\eta}}(\eta) \rightarrow 1 \ \ \mbox{as}
\ \ \eta \rightarrow \infty $.
For this particular value of $ \lambda $, due to the scaling properties, we have found that
\begin{equation}
f_M(\eta) = \lambda f(\lambda \eta)
\end{equation}
because $ \lambda f(\lambda \eta) $ verifies the BVP (\ref{p7}) that defines $ f_M(\eta) $ uniquely.

Hence, the error for $ \eta \in [0, M] $ is given by
\begin{equation}
e(\eta) = |\lambda f(\lambda \eta) - f(\eta)| \leq
|(\lambda -1) f(\lambda \eta)| + |f(\lambda \eta) - f(\eta)| \quad .
\end{equation}
By applying the mean value theorem of differential calculus and taking
into account that $ {\displaystyle \frac{df}{d\eta}}(\eta) \leq 1 $
we get the relations $ f(\lambda \eta) \leq \lambda \eta $ and
$ |f(\lambda \eta) - f(\eta)| \leq (\lambda - 1) \eta $. As a result
\begin{equation}
e(\eta) \leq M (\lambda^2 -1), \qquad \eta \in [0, M]
\end{equation}
where $ \lambda^2 -1 > 0 $ because $ \lambda > 1 $.
Naturally,
$ {\displaystyle \frac{df_M}{d\eta}}(\eta \rightarrow \infty) =
\lambda^2 $, so that
\begin{eqnarray*}
& \begin{array}{ll}
\lambda^2 -1 &= {\displaystyle \frac{df_M}{d\eta}}(\eta \rightarrow \infty)
- {\displaystyle \frac{df_M}{d\eta}}(M) \\
& = {\displaystyle \int_{M}^{\infty}}
{\displaystyle \frac{d^2 f_M}
{d\eta^2}} d\eta \quad .
\end{array}
\end{eqnarray*}
To complete the proof Rubel used some manipulations,
involving a first integral of the governing differential equation, to find that
\begin{equation}
\lambda^2 -1 \leq {\displaystyle \frac{d^2 f_M}{d\eta^2}}(M)
[f_M(M)]^{-1} \ .
\end{equation}
\hfill $ \Box $

\noindent
{\bf Remark.} As a consequence of this theorem in order to control the error we can modify either the value of $ M $ or the value of $ {\displaystyle \frac{d^2 f_M}{d\eta^2}}(M) $. 
Classically the value of $ M $ has been chosen to this end.
The above Theorem shows that the error is directly proportional to $ M $. 
In this context Fazio defines a free boundary formulation of the Blasius problem where the second order derivative of the solution with respect to $ \eta $ at the free boundary can be chosen as small as possible, see \cite{Fazio:1992:BPF} for details.
The free boundary can be interpreted as a truncated boundary, see Fazio \cite{Fazio:2002:SFB}.

%% file: Tab_Blasius_Slip_Iter.tex
\begin{table}[!htb]
\renewcommand\arraystretch{1.5}
	\centering
		\begin{tabular}{r@{.}lr@{.}lr@{.}lr@{.}l}
\hline
\multicolumn{2}{c}%
{$ c $}
& \multicolumn{2}{c}%
{${\displaystyle \frac{df^*}{d\eta^*}(\eta^*_{\infty})}$}
& \multicolumn{2}{c}%
{${\displaystyle \frac{df}{d\eta}(0)}$}
& \multicolumn{2}{c}%
{${\displaystyle \frac{d^2f}{d\eta^2}(0)}$} \\ \hline
 0  & & 2  & 085393 & 0 &        & 0 & 332061    \\
 1  &  & 2  & 262516 & 0 & 293841 & 0 & 293841    \\
 5  & & 3  & 644351 & 0 & 718686 & 0 & 143737    \\
10  &  & 5  & 203210 & 0 & 842545 & 0 & 084255    \\
50  &  &13  & 894469 & 0 & 965399 & 0 & 019308    \\
\hline			
		\end{tabular}
	\caption{Slip boundary condition: numerical results by the ITM.}
	\label{tab:Piter}
\end{table}

%% file: Tab_Blasius_Move_Iter.tex
\begin{table}[!htb]
\renewcommand\arraystretch{1.5}
	\centering
		\begin{tabular}{r@{.}lr@{.}lr@{.}lr@{.}l}
\hline 
\multicolumn{2}{c}%
{$h^*$}
& \multicolumn{2}{c}%
{$\Gamma(h^*)$}
& \multicolumn{2}{c}%
{$h^*$} 
& \multicolumn{2}{c}%
{$ \Gamma(h^*) $}\\ 
\hline 
 1 &        & $-$0 & 528765 &          10 &         & $-$0 & 008382 \\
 3 &        &  0 & 198514 &           3 &         &  0 & 198514 \\
 2 & 454091 &  0 & 045354 &           9 & 716410  &  0 & 016266 \\
 2 & 339221 &  0 & 008091 &           9 & 903562  & $-$1 & $52\mbox{D-04}$ \\
 2 & 319037 &  0 & 001371 &           9 & 901831  & $-$2 & $96\mbox{D-06}$ \\
 2 & 315626 &  2 & $21\cdot 10^{-4}$ & 9 & 901798  & $-$5 & $79\mbox{D-08}$ \\
 2 & 315076 &  3 & $87 \mbox{D-05}$  &  &  &  &  \\
 2 & 314979 &  6 & $78 \mbox{D-06}$  &  &  &  &  \\
 2 & 314963 &  1 & $19 \mbox{D-06}$ &  &  &  &  \\
 2 & 314960 &  2 & $08 \mbox{D-07}$  &  &  &  &  \\
\hline
\multicolumn{2}{c}%
{${\displaystyle \frac{df}{d\eta}(0)}$}
& \multicolumn{2}{c}%
{${\displaystyle \frac{d^2f}{d\eta^2}(0)}$}
& \multicolumn{2}{c}%
{${\displaystyle \frac{df}{d\eta}(0)}$} 
& \multicolumn{2}{c}%
{$ {\displaystyle \frac{d^2f}{d\eta^2}(0)}$}\\ \hline
 $-$0 & 25 & 0 & 283928  & $-$0 & 25 & 0 & 032094 \\
\hline			
		\end{tabular}
	\caption{Fluid flow on a moving surface for $b = - 0.25$: numerical results by the ITM.}
	\label{tab:Move}
\end{table}

%% file: FS.tex
\section{Fluid flow on wedges}
\begin{figure}[!hbt]
	\centering
\begin{small}
\psfrag{y}[][]{$ y $} 
\psfrag{z}[][]{$ z $} 
\psfrag{O}[][]{$ O $}
\psfrag{V}[][]{$ V_\infty $} 
\includegraphics[width=.48\textwidth]{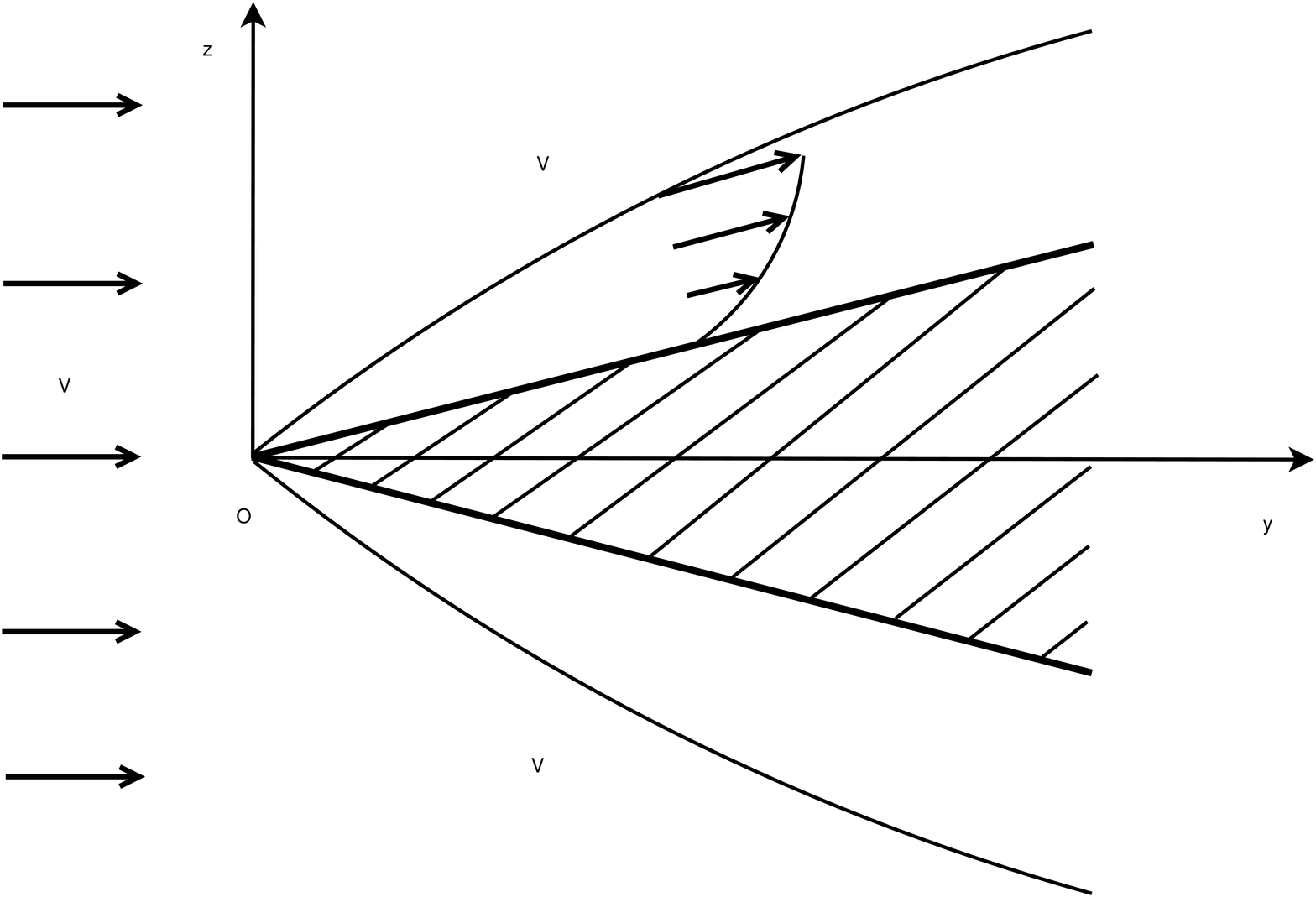}
\hfil
\includegraphics[width=.48\textwidth]{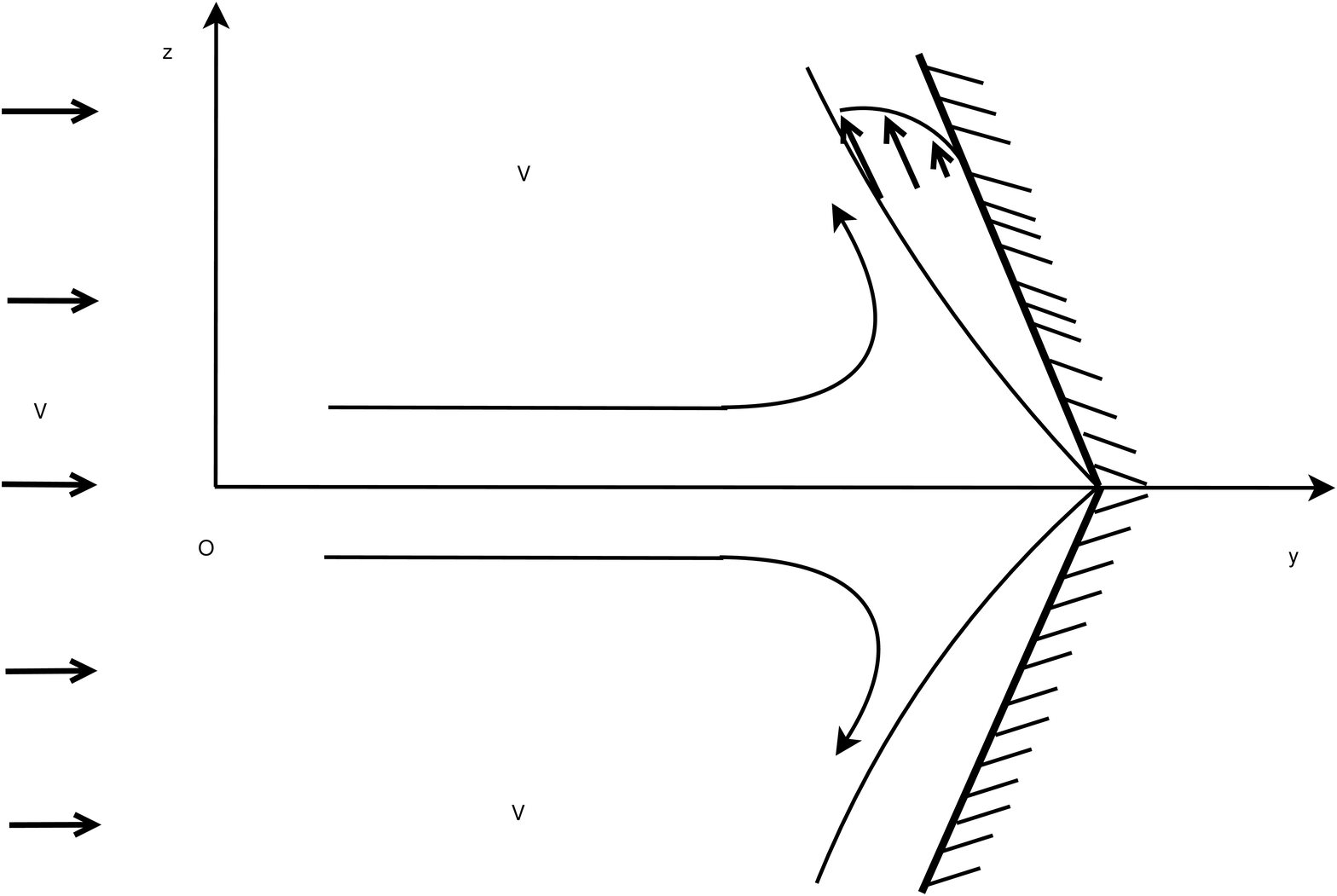}
\hfil
\end{small}
\caption{Boundary layer over different wedges.}
	\label{fig:Wedges}
\end{figure}

Within the celebrated boundary-layer theory, developed at the beginning of the last century by Prandtl \cite{Prandtl:1904:UFK}, the model describing the steady plane flow of a fluid past a wedge, provided the boundary layer assumptions are verified ($ v \gg w $ and the existence of a very thin layer attached to the wedge), is given by
\begin{align}
& {\displaystyle \frac{\partial v}{\partial y}} +
{\displaystyle \frac{\partial w}{\partial z}} = 0  \\
& v {\displaystyle \frac{\partial v}{\partial y}} +
w {\displaystyle \frac{\partial v}{\partial z}} =
V_{\infty} {\displaystyle \frac{dV_{\infty}}{dy}} + \mu
{\displaystyle \frac{\partial^2 v}{\partial z^2}}  \\
& v(y, 0) = w(y, 0) = 0 \\
& v(y, z) \rightarrow V_{\infty}(y)
\quad \mbox{as} \quad z \rightarrow \infty
\end{align}
where the governing differential equations, conservation of mass and momentum, are the steady-state 2D Navier-Stokes equations under the boundary layer approximations, $ v $ and $ w $ are the velocity components of the fluid in the $ y $ and $ z $ direction, $ V_{\infty} (y) $ represents the main-stream velocity, see the draft in figure \ref{fig:Wedges}, and $ \mu $ is the viscosity of the fluid.
The boundary conditions at $ z = 0 $ are based on the assumption that neither slip nor mass transfer is permitted, whereas the remaining boundary condition means that the velocity $ v $ tends to the main-stream velocity $ V_{\infty}(y) $ asymptotically.

In order to study this problem, it is convenient to introduce a potential (stream function) $ \psi(y, z) $ defined by
\begin{equation}
v = {\displaystyle \frac{\partial \psi}{\partial z}}, \qquad
w = - {\displaystyle \frac{\partial \psi}{\partial y}} \quad .
\end{equation}
After that the problem becomes
\begin{align}
& \mu {\displaystyle \frac{\partial^3 \psi}{\partial z^3}} +
V_{\infty} {\displaystyle \frac{dV_{\infty}}{dy}} +
{\displaystyle \frac{\partial \psi}{\partial y}}
{\displaystyle \frac{\partial^2 \psi}{\partial z^2}} 
-{\displaystyle \frac{\partial \psi}{\partial z}}
{\displaystyle \frac{\partial^2 \psi}{\partial y \partial z}}
= 0  \\
& {\displaystyle \frac{\partial \psi}{\partial y}}(y, 0) =
{\displaystyle \frac{\partial \psi}{\partial z}} (y, 0) = 0  \\
& {\displaystyle \frac{\partial \psi}{\partial z}} (y, z)
\rightarrow V_{\infty} (y) \quad \mbox{as} \quad
z \rightarrow \infty \quad .
\end{align}
A similarity analysis shows that it must be $ V_{\infty} (y) \propto y^m $ and that (taking the constant of proportionality equal to one) by introducing the new variables 
\[
 \eta = \left[ \frac{1}{2 \mu} (1 + m) \right]^{1/2} z y^{-(1-m)/2} 
\]
and 
\[
 f(\eta) = \left[ \frac{1}{2 \mu} (1 + m) \right]^{1/2}
\psi(y, z) y^{-(m+1)/2} 
\]
our model reduces to
\begin{eqnarray}\label{eq:Falkner}
& {\displaystyle \frac{d^3 f}{d \eta^3}} + f
{\displaystyle \frac{d^{2}f}{d\eta^2}} + \beta \left[ 1 -
\left({\displaystyle
\frac{df}{d\eta}}\right)^2 \right] = 0 \nonumber \\[-1.5ex]
& \\[-1.5ex]
& f(0) = {\displaystyle \frac{df}{d\eta}}(0) = 0, \qquad
{\displaystyle \frac{df}{d\eta}}(\eta) \rightarrow 1 \quad \mbox{as}
\quad \eta \rightarrow \infty \nonumber
\end{eqnarray}
where $ \beta = 2m/(1 + m) $.
The governing differential equation in (\ref{eq:Falkner}) is known as the Falkner-Skan equation \cite{Falkner:1931:SAS}.
It has been proved by Weyl \cite{Weyl:1942:DES} that for each value of $ \beta $ there exists a solution for which its second derivative is positive, monotone decreasing on $ [0, \infty) $ and approaching zero as $ \eta $ goes to infinity.
The uniqueness question is more complex: if $ \beta > 1 $, then besides the monotone solution a hierarchy of solution with reversed flow exists as shown by Coppel \cite{Coppel:1960:DEB}, see also Craven and Pelietier \cite{Craven:1972:RFS}.
The classical case, corresponding to $ \beta = m = 0 $, is known as the Blasius problem \cite{Blasius:1908:GFK}.

\section{The Falkner-Skan model}
In this section, we apply the ITM to the Falkner-Skan equation with relevant boundary conditions
\begin{align}\label{eq:abf1}
& {\displaystyle \frac{d^{3}f}{d\eta^3}} + f 
{\displaystyle \frac{d^{2}f}{d\eta^2}} + \beta \left[ 1 - \left({\displaystyle
\frac{df}{d\eta}}\right)^2 \right] = 0 \nonumber \\[-1.2ex]
& \\[-1.2ex]
& f(0) = {\displaystyle \frac{df}{d\eta}}(0) = 0 \ , \qquad
{\displaystyle \frac{df}{d\eta}}(\eta) \rightarrow 1 \quad \mbox{as}
\quad \eta \rightarrow \infty \ , \nonumber
\end{align}
where $ f $ and $ \eta $ are appropriate similarity variables
and $ \beta $ is a parameter. 
This set is called the Falkner-Skan model, after the names of two English mathematicians
who first studied it \cite{Falkner:1931:SAS}.
As pointed out by Na \cite[pp. 146-147]{Na:1979:CME}, if $ \beta \ne 0$, the BVP (\ref{eq:abf1}) cannot be solved by a non-ITM.
Indeed, the governing differential equation in (\ref{eq:abf1}) is not invariant with respect to any scaling group of point transformations.

The existence and uniqueness question for the problem (\ref{eq:abf1}) is really a complex matter.
Assuming that $ \beta > 0 $ and under the restriction $ 0 < \frac{df}{d\eta} < 1 $, known as normal flow condition, Hartree \cite{Hartree:1937:EOF} and Stewartson \cite{Stewartson:1954:FSF} proved that the problem (\ref{eq:abf1}) has a unique solution, whose first derivative tends to one exponentially.
Coppel \cite{Coppel:1960:DEB} and Craven and Peletier \cite{Craven:1972:USF} pointed out that the above restriction on the first derivative can be omitted when $ 0 \le \beta \le 1$.
Weyl proved, in \cite{Weyl:1942:DES}, that for each value of the parameter $\beta$ there exists a physical solution with positive  monotone decreasing, in $[0, \infty)$, second derivative that approaches zero as the independent variable goes to infinity. 
In the case $ \beta > 1 $, the Falkner-Skan model loses the uniqueness and a hierarchy of solutions with reverse flow exists.
In fact, for $ \beta > 1$ Craven and Peletier \cite{Craven:1972:RFS} computed solutions for which $ \frac{df}{d\eta} < 0 $ for some value of $ \eta $.
In each of these solutions the velocity approaches its limit exponentially in $\eta$. 
As mentioned before, the term normal flow indicates that
the flow velocity has a unique direction, and instead, reverse flow means that the velocity
is both positive and negative in the integration interval. 

The considered problem, has also multiple solutions for $\beta_{\min} < \beta < 0$, as reported by  Veldman and van de Vooren \cite{Veldman:1980:GFS}, with the minimum value of $\beta$ given by 
\begin{equation}\label{eq:bmin}
\beta_{\min} = -0.1988\dots \ .
\end{equation} 
In this range, there exist two physical solutions, one for normal flow and one
for reverse flow.  
For $\beta = \beta_{\min}$ only one solution exists.
Finally, for $\beta < \beta_{\min}$ the problem has no solution at all. 
Our interest here is to apply the ITM to the range of $\beta$ where multiple solutions are admitted, in particular, our interest is to get numerically the famous solutions of Stewartson \cite{Stewartson:1954:FSF,Stewartson:1964:TLB}.
The obtained results are original and, as we shall see in the following sections, are in agreement with those available in the literature.

The first computational treatment of the Falkner-Skan model is due to Hartree \cite{Hartree:1937:EOF}.
Cebeci and Keller \cite{Cebeci:1971:SPS} apply shooting and parallel shooting methods requiring asymptotic boundary condition to be imposed at a changing unknown boundary in the computation process. 
As a result, they report convergence difficulties, which can be avoided by moving towards more complicated methods. 
Moreover, to guarantee reasonable accuracy, they are forced to use a {\it small enough} step-size and extensive computation for the solution of the IVPs.
Na \cite[pp. 280-286]{Na:1979:CME} describes the application of invariant imbedding.
A modified shooting method \cite{Asaithambi:1997:NMS} and finite-difference methods \cite{Asaithambi:1998:FDM,Asaithambi:2004:SOF} for this problem are presented by Asaithambi.
Kuo \cite{Kuo:2003:ADT} uses a differential transformation method, which obtains a series solution of the Falkner-Skan equation.
Sher and Yakhot \cite{Sher:2001:NAS} define a new approach that solves this problem by shooting from infinity, using some simple analysis of the asymptotic behaviour of the solution at infinity.
Asaithambi \cite{Asaithambi:2005:SFE} proposes a faster shooting method by using recursive evaluation of Taylor coefficients.
Zhang and Chen \cite{Zhang:2009:IMS} investigate a modification of the shooting method, where the computation of the Jacobian matrix is obtained by solving two IVPs. 
A Galerkin-Laguerre spectral method is defined by Auteri and Quartapelle \cite{Auteri:2012:GLS}.

\subsection{The ITM}
In order to apply an ITM to (\ref{eq:abf1}) we have to embed it in a modified model and require the invariance of this last model with respect to an extended scaling group of transformations.
This can be done in several ways that are all equivalent.
In fact, the modified model can be written as 
\begin{align}\label{eq:abf2}
& {\displaystyle \frac{d^{3}f}{d\eta^3}} + f 
{\displaystyle \frac{d^{2}f}{d\eta^2}} + \beta \left[ h^{4/\sigma} - \left({\displaystyle
\frac{df}{d\eta}}\right)^2 \right] = 0 \ , \nonumber \\[-1.2ex]
& \\[-1.2ex]
& f(0) = {\displaystyle \frac{df}{d\eta}}(0) = 0 \ , \qquad
{\displaystyle \frac{df}{d\eta}}(\eta) \rightarrow 1 \quad \mbox{as}
\quad \eta \rightarrow \infty \ , \nonumber
\end{align}
and the related extended scaling group is given by
\begin{equation}\label{eq:scalinv}
f^* = \lambda f \ , \qquad \eta^* = \lambda^{-1} \eta \ , \qquad 
h^* = \lambda^{\sigma} h \ ,   
\end{equation}
where $\sigma$ is a parameter.
In the following we set $\sigma = 4$; for the choice $\sigma=8$ see \cite{Fazio:1994:FSE}.
In \cite{Fazio:1994:FSE}, a free boundary formulation of the Falkner-Skan model was considered and
numerical results were computed for the Homann flow ($ \beta = 1/2 $) as well as for the Hiemenz flow ($ \beta = 1 $).

From a numerical point of view the request to evaluate 
$ \frac{d f}{d \eta} (\infty) $ cannot be fulfilled.
Several strategies have been proposed in order to provide an approximation of this value.
The simplest and widely used one is to introduce, instead of infinity, a suitably truncated boundary.
A recent successful way to deal with such an issue is to reformulate the considered problem as a free BVP \cite{Fazio:1992:BPF,Fazio:1994:FSE,Fazio:1996:NAN}; for a survey on this topic see \cite{Fazio:2002:SFB}. 
Recently, Zhang and Chen \cite{Zhang:2009:IMS} have used a free boundary formulation to compute the normal flow solutions of the Falkner-Skan model in the full range $ \beta_{\min} < \beta \le 40$.
They applied a modified Newton's method to compute both the initial velocity and the free boundary.  
For the sake of simplicity, we do not use the free boundary approach but,  following T\"opfer,  we use some preliminary computational tests to find a suitable value for the truncated boundary.

At each iteration of the ITM, we have to solve the IVP 
\begin{align}\label{eq:IVP2}
& {\displaystyle \frac{d^{3}f^*}{d\eta^{*3}}} + f^* 
{\displaystyle \frac{d^{2}f^*}{d\eta^{*2}}} + \beta \left[ h_j^{*} - \left({\displaystyle
\frac{df^*}{d\eta^*}}\right)^2 \right] = 0 \nonumber \\[-1.2ex]
& \\[-1.2ex]
& f^*(0) = {\displaystyle \frac{df^*}{d\eta^*}}(0) = 0 \ , \qquad
{\displaystyle \frac{d^2f^*}{d\eta^{*2}}}(0) = \pm 1  \ . \nonumber
\end{align}
Tables~\ref{tab:Itera1} and~\ref{tab:Itera2} list the numerical iterations obtained for a sample value of $\beta$.
We notice that we solve an IVP governed by a different differential equation for each iteration because the Falkner-Skan equation is not invariant under every scaling group of point transformation.
We have chosen $ \beta=-0.01$ since, in this case, the missing initial conditions for the normal and reverse flows are not symmetric with respect to the $\beta$ axis. 

\input{CF2013TAB1}

The data listed in tables~\ref{tab:Itera1} and~\ref{tab:Itera2} have been obtained by solving the modified Falkner-Skan model on $ \eta^* \in [0, 20] $ by setting 
\[
 \frac{d^2f^*}{d\eta^{*2}}(0) = \pm 1 \ ,
\]
respectively.

\input{CF2013TAB2}

In both cases, we achieved convergence of the numerical results within seven iterations.
Let us now investigate the behaviour of the transformation function.
\begin{figure}[!hbt]
	\centering
\psfrag{h*}[][]{$h^*$} 
\psfrag{G}[][]{$\Gamma(h^*)$} 
\includegraphics[width=.475\textwidth]{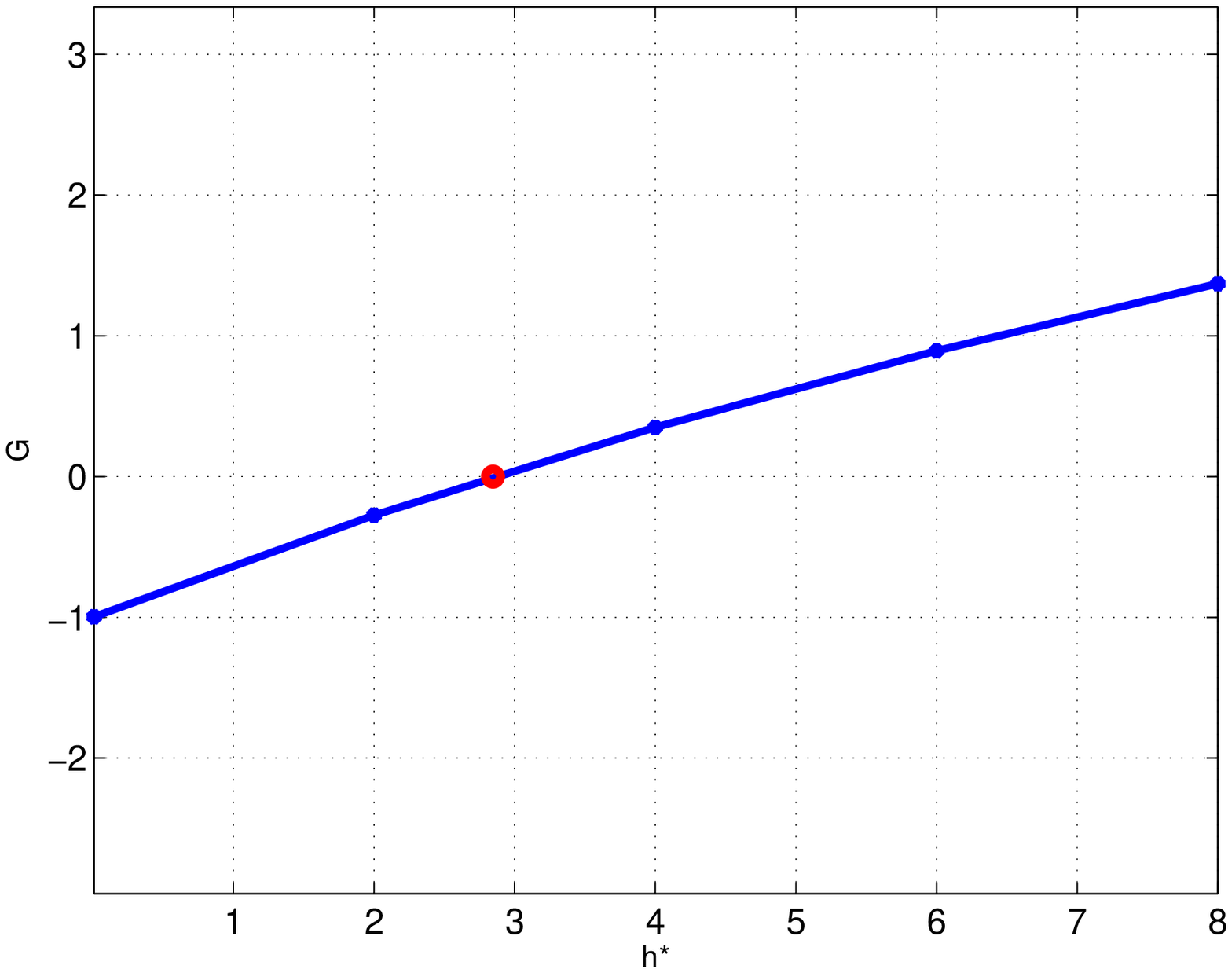}
\hfill
\includegraphics[width=.475\textwidth]{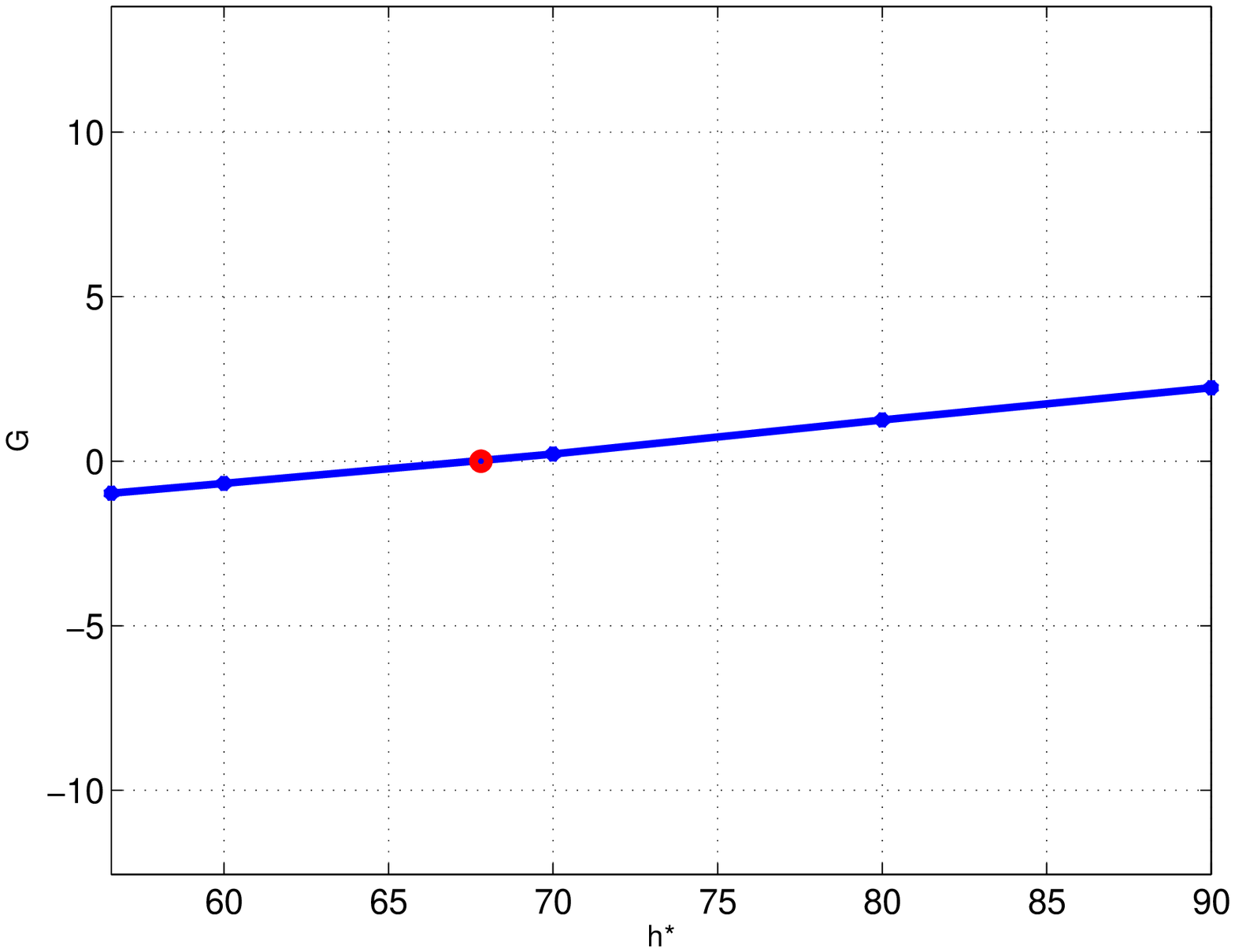}
\caption{Two cases of the $\Gamma(h^*)$ function: left and right frames are related to normal and reverse flow solutions, respectively.} 
	\label{fig:FSHGhstar}
\end{figure}
Figure~\ref{fig:FSHGhstar} shows $\Gamma(h^*)$ with respect to $h^*$ for the two cases reported in these tables.
The unique zero of the transformation function is marked by a circle.
It is worth noticing that the same scale has been used for both axes.
As it is easily seen, in both cases, we have a monotone increasing function.
We notice on the left frame, corresponding to a normal flow, that the tangent to the $\Gamma$ function at its unique zero and the $h^*$ axis, define a large angle.
This is important from a numerical viewpoint because in such a case we face a well-conditioned problem.
On the other hand, this is not the case for the function plotted on the right frame of the same figure. 
The meaning is clear, reverse flow solutions are more challenging to compute than normal flow ones.
Therefore, one has to put some care when choosing the convergence criteria for the root-finder method.
\begin{figure}[!hbt]
	\centering
\psfrag{e}[][]{\small $ \eta $} 
\psfrag{f}[][]{\small $ f(\eta)$} 
\psfrag{fe}[][]{$ {\displaystyle \frac{df}{d\eta}} $} 
\psfrag{fee}[][]{$ {\displaystyle \frac{d^2f}{d\eta^2}} $} 
\includegraphics[width=\tw,height=6cm]{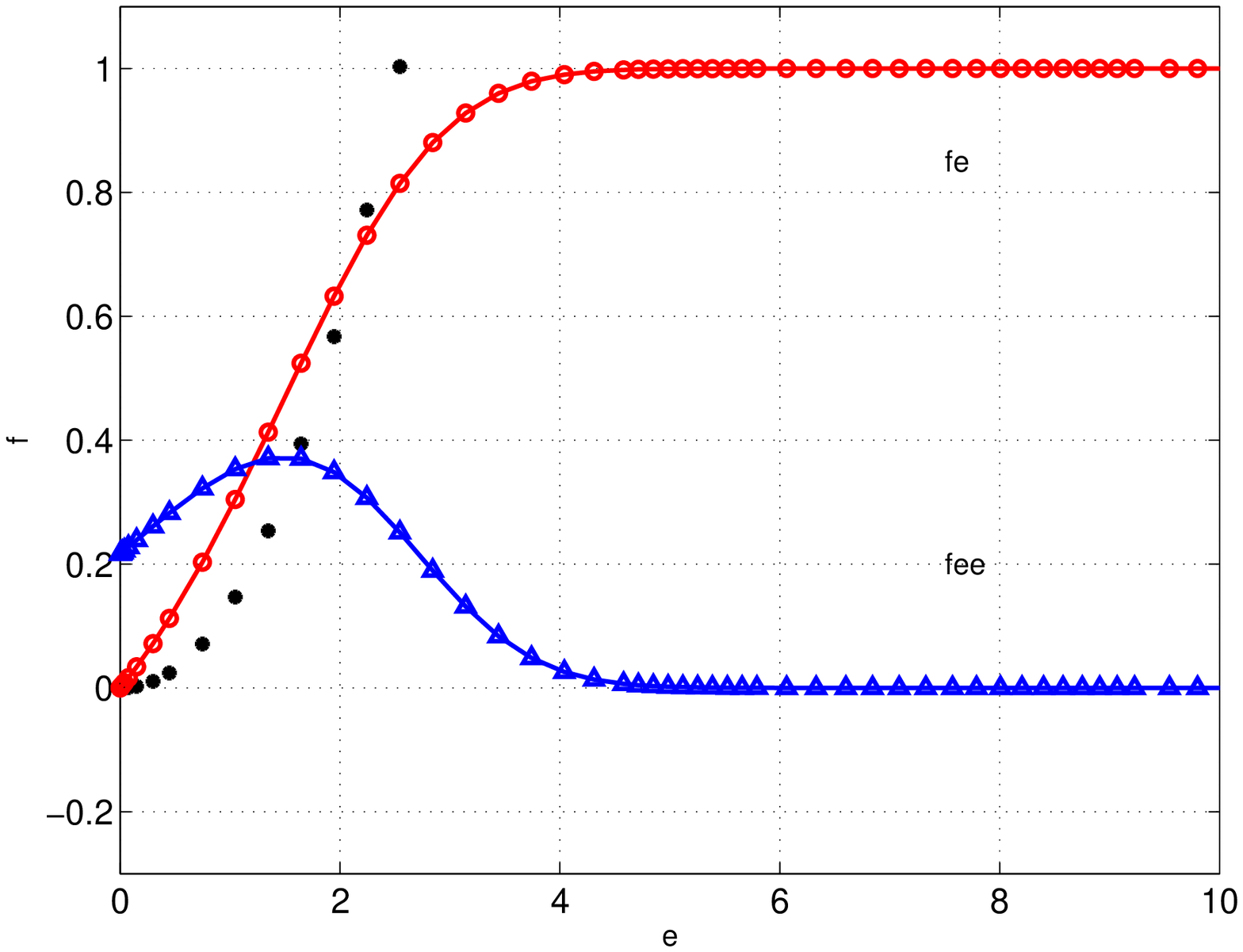} \\
\includegraphics[width=\tw,height=6cm]{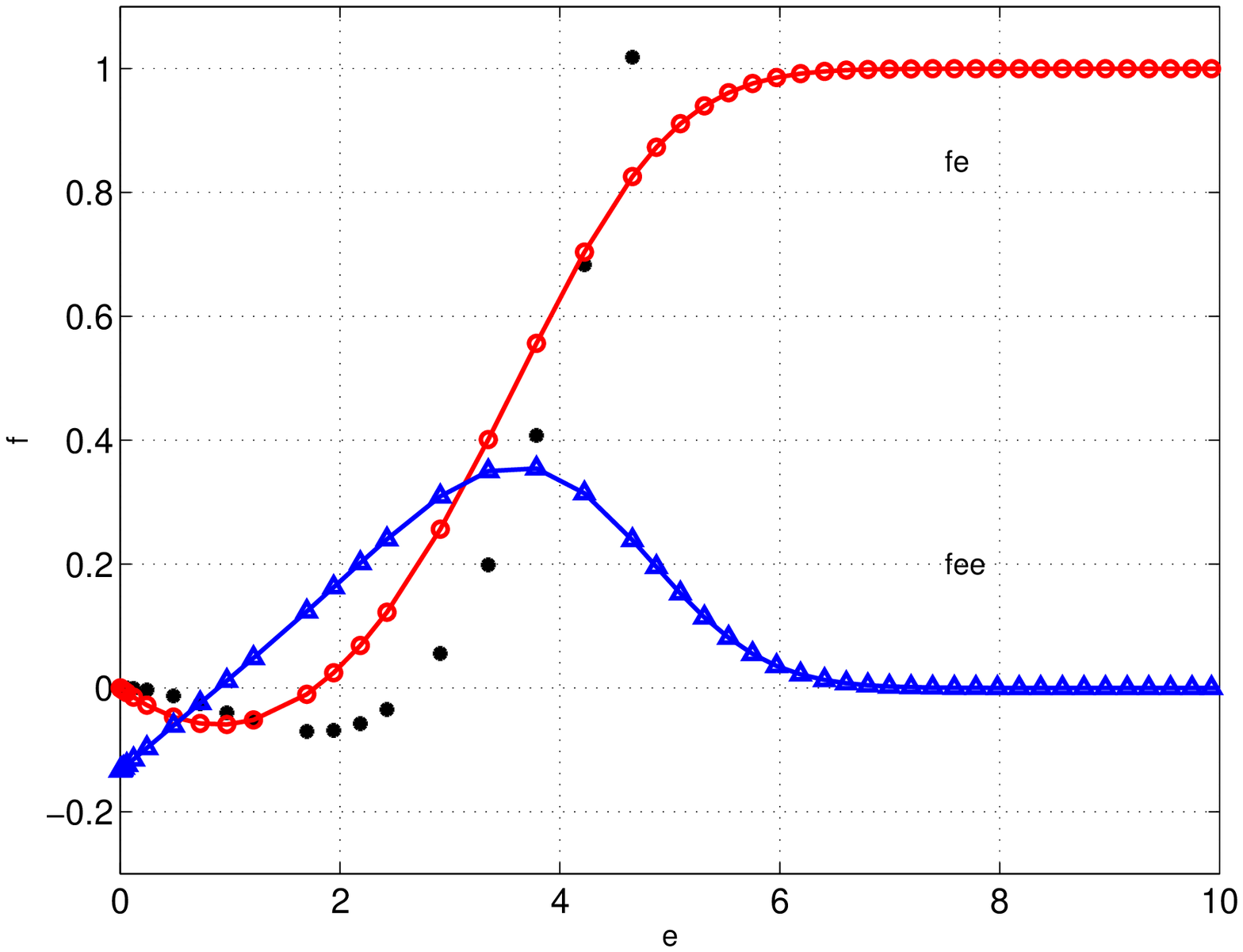}
\caption{Normal and reverse flow solutions to Falkner-Skan model for $\beta =-1.5$.
The symbols $\bullet$ denote values of $f(\eta)$.} 
	\label{fig:Falkner}
\end{figure}
Figure~\ref{fig:Falkner} shows the results of the two numerical solutions for a different value 
of $\beta$, namely $\beta = -0.15$.
In the top frame, we have the normal flow and in the bottom frame we display the reverse flow solution.
In both cases the solutions were computed by introducing a truncated boundary and solving the IVP in the starred variables on $\eta^* \in [0, 20]$ with $h_0^*=1$, $h_1^* = 5$ in the top frame and
$h_0^*=15$, $h_1^* = 25$ in the bottom frame.
In this case, we achieved convergence of the numerical results within eight and seven iterations, respectively.
For the sake of clarity, we omit to plot the solutions in the starred variables computed during the iterations.
Moreover, we display only $\eta \in [0, 10]$.


\input{CF2013TAB4}
As far as the reverse flow solutions are concerned, in table~\ref{tab:Compare} we compare the missing initial condition computed by the ITM for several values of $\beta$ with results available in the literature.
The agreement is really good.
It is remarkable that among the studies quoted in the introduction only a few report data related to the reverse flow solutions. 

\begin{figure}[!hbt]
	\centering
\psfrag{b}[][]{$ \beta $} 
\psfrag{Blasius}[l][l]{Blasius} 
\psfrag{fee}[][]{${\displaystyle \frac{d^2f}{d\eta^2}(0)} $} 
\includegraphics[width=\tw]{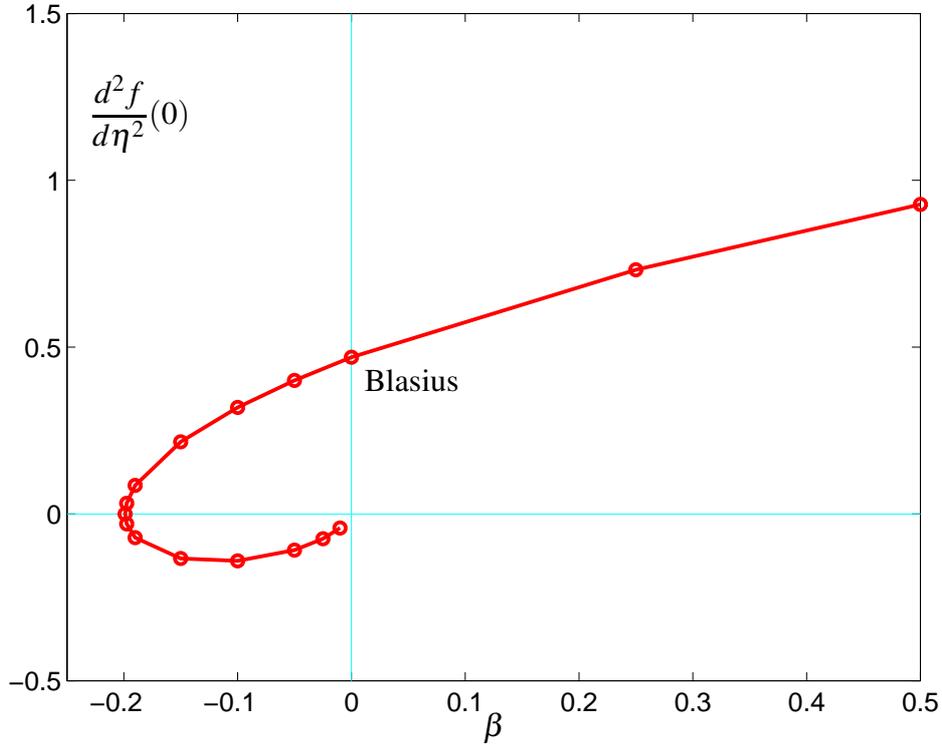}
\caption{Missing initial conditions to Falkner-Skan model for several values of $\beta$.
Positive values determine normal flow, and instead, negative values define reverse flow solutions.} 
	\label{fig:bfee}
\end{figure}
In figure~\ref{fig:bfee} we plot the behaviour of missed initial condition versus $\beta$.
The solution found by the data in table~\ref{tab:Itera2} is plotted in this figure, but not the one found in table~\ref{tab:Itera1} because this is very close to the Blasius solution.
A good initial choice of the initial iterates of $h^*$, for a given value of $\beta$, is obtained by employing values close to the one used in a successful attempt made for a close value of $\beta$.
It is interesting to note that, for values of $\beta < \beta_{\min}$ the ITM continued to iterate endlessly, whatever set of starting values for $h^*$ are selected.

Our extended algorithm has shown a kind of robustness because it is able to get convergence even when, for a chosen value of $h^*$, the IVP solver stops before arriving at the selected truncated boundary getting a wrong value of $\Gamma(h^*) = -1$.
On the other hand, the secant method gives an overflow error when this happens for two successive iterates of $h^*$.

The value of $\beta_{\min}$, corresponding to a separation point at $\eta=0$, can be found by the ITM by considering $\beta$ as a continuation parameter.
As we have seen, for $ \beta_{\min} < \beta < 0$ two solutions are available: a positive and a negative skin-friction coefficient, the missing initial condition, providing a normal and reverse flow solution.
For instance, when $\beta=-0.1988$ we get for the missing initial conditions the values
$0.005221$ and $-0.005158$, respectively.
Starting from this value of $\beta$ we can reduce it gradually and check whether the two missing initial conditions, the positive and negative values of $ {\frac{d^2f}{d\eta^2}(0)}$, converge to zero.  
Soon, we realize that we are forced the use the ITM to its natural limit.
In fact, we are trying to get a skin-friction coefficient close to zero rescaling a fixed non-zero value, plus or minus one in our case.
Anyway, when $\beta = -0.198837723795$ we found the skin-friction coefficients $6.61\mbox{D}-06$ and $-6.61\mbox{D}-06$ with 20 and 24 iterations, respectively.
Finally, we have noticed that, as far the guest for this limiting value of $\beta$ is concerned, we are allowed to reduce the chosen truncated boundary value, and for $\beta = -0.198837723795$ this truncated boundary was set equal to one, i.e. all IVP was solved on $[0, 1]$.

\begin{figure}[!hbt]
	\centering
\psfrag{e}[][]{\small $ \eta $} 
\psfrag{f}[][]{\small $f(\eta)$} 
\psfrag{fe}[][]{$ {\displaystyle \frac{df}{d\eta}} $} 
\psfrag{fee}[][]{$ {\displaystyle \frac{d^2f}{d\eta^2}} $} 
\includegraphics[width=\tw]{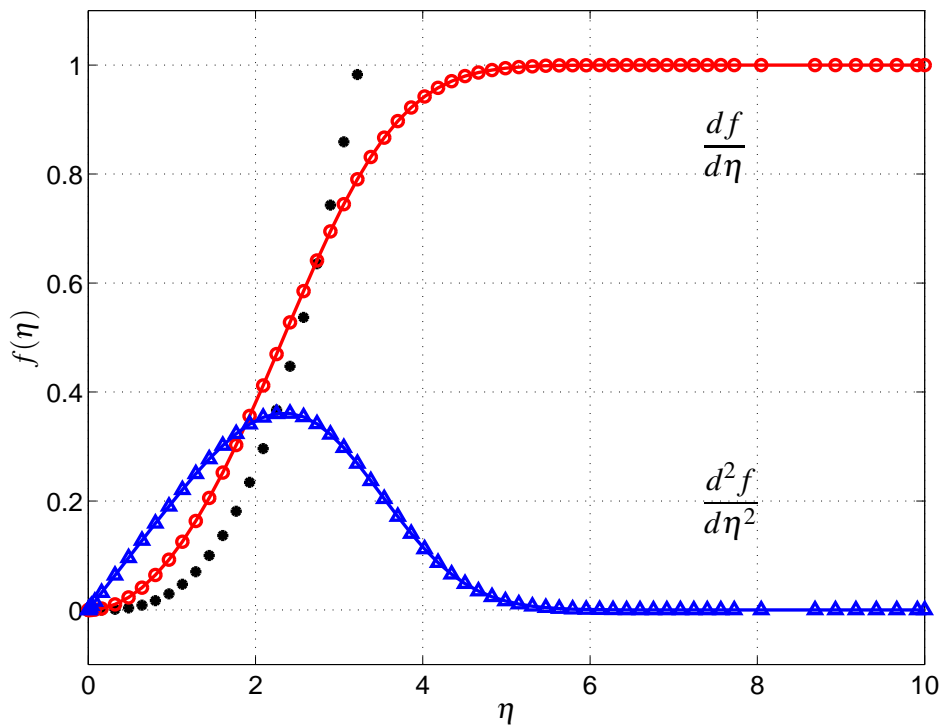}
\caption{Numerical solutions to Falkner-Skan model for $\beta = -0.1988376$. We notice that $ {\displaystyle \frac{d^2f}{d\eta^2}}(0) = 0$ and values of $f(\eta)$ are marked by $\bullet$.} 
	\label{fig:FSHL}
\end{figure}
In figure~\ref{fig:FSHL} we plot the unique solution for the limiting value $\beta_{\min}$, where $\beta_{\min}$ is given by equation (\ref{eq:bmin}).
As it easily seen this is a normal flow solution. 

\begin{figure}[!htb]
	\centering
\psfrag{e}[][]{\small $ \eta $} 
\psfrag{fe}[l][]{$ {\frac{df}{d\eta}} $} 
\psfrag{fee}[][l]{$ {\frac{d^2f}{d\eta^2}} $} 
\includegraphics[width=.45\textwidth,height=\textwidth]{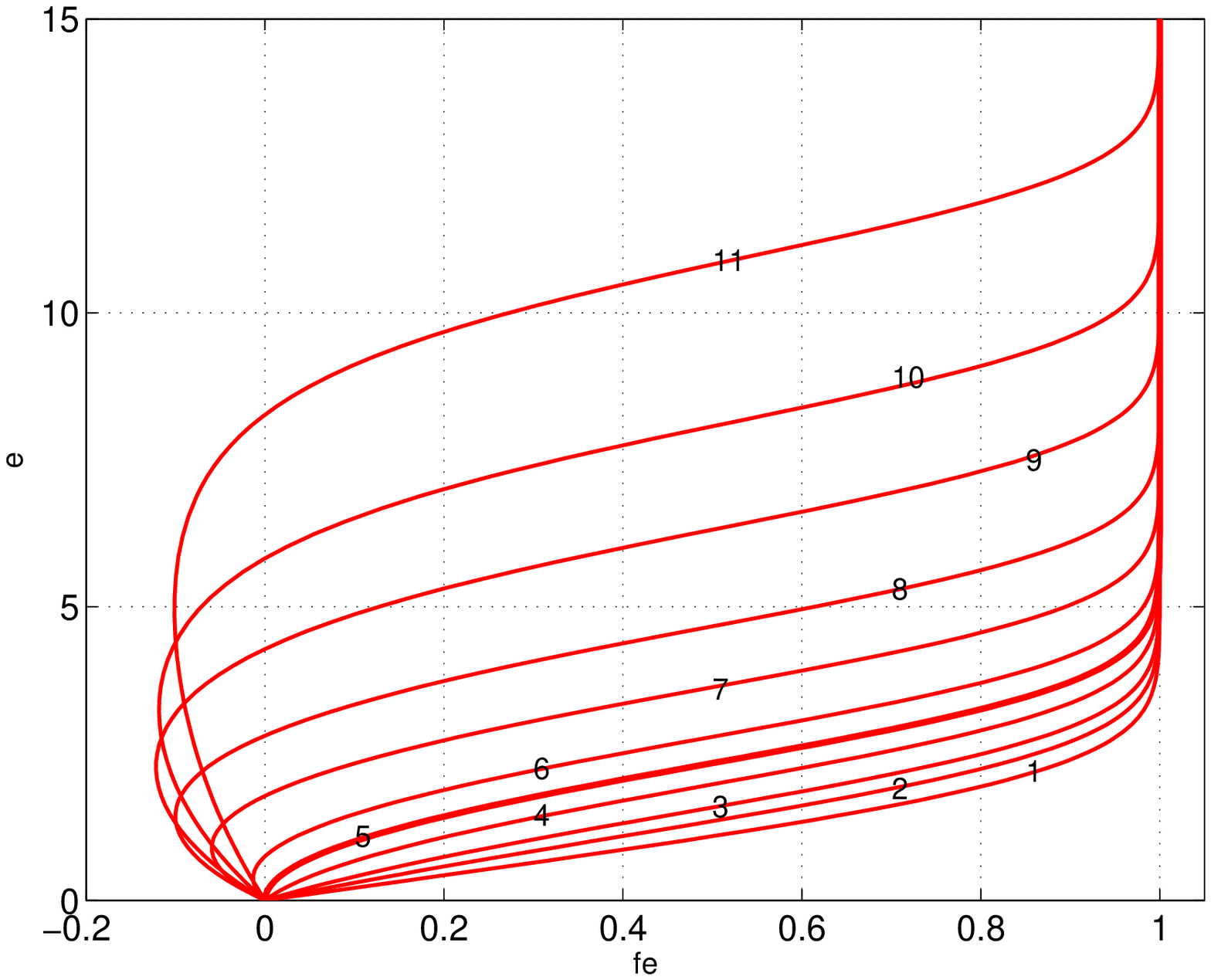} \hfill
\includegraphics[width=.45\textwidth,height=\textwidth]{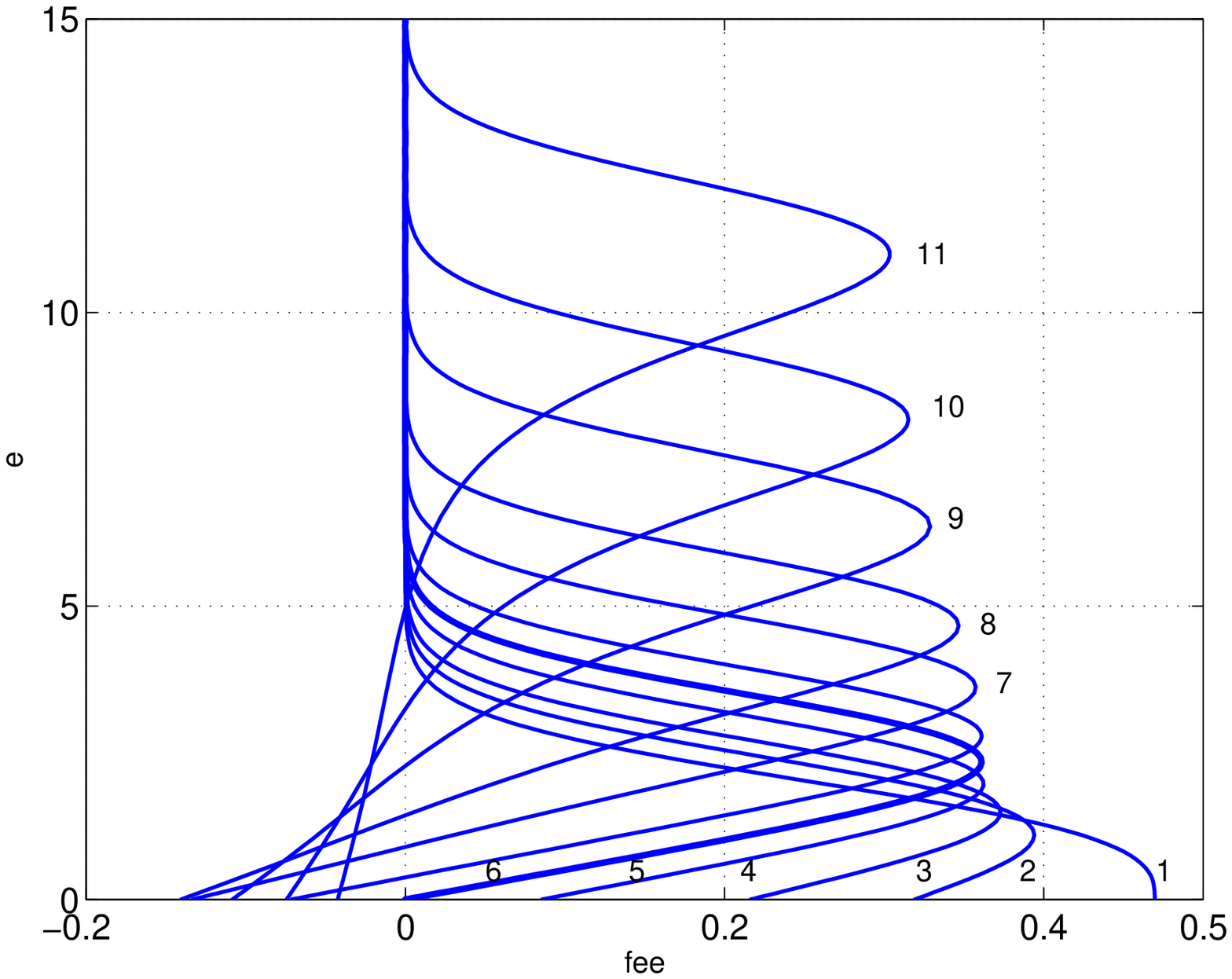}
\caption{Numerical solutions to Falkner-Skan model. On the left frame we plot $ {\frac{df}{d\eta}}(\eta)$ and on the right frame we show $ {\frac{d^2f}{d\eta^2}}(\eta)$.
A different number indicate the corresponding value of $\beta$: 1 stands for $\beta = 0$, the Blasius solution, 2 for $\beta = -0.1$, 3 for $\beta = -0.15$, 4 for $\beta = -0.19$, 5 for $\beta_{\min}$, 6 for $\beta = -0.19$, 7 for $\beta = -0.15$, 8 for $\beta = -0.1$, 9 for $\beta = -0.05$, 10 for $\beta = -0.025$, 11 for $\beta = -0.01$.} 
	\label{fig:FSHplot}
\end{figure}

Figure \ref{fig:FSHplot} shows the numerical results for several values of $\beta$.
The solution corresponding to the limiting value $\beta_{\min}$ is marked by a heavy line.

The results reported so far have been found by a variable order adaptive multi-step IVP solver that was coupled up the simple secant method.
The adaptive solver uses a relative and an absolute error tolerance, for each component of the numerical solution, both equal to ten to the minus six. 
As well known, the secant method is convergent provided that two initial iterates sufficiently close to the root are used, and its convergence is superlinear with an order of convergence equal to $(1+\sqrt{5})/2$.
As far as a termination criterion for the secant method is concerned, we enforced the conditions
\begin{equation}\label{eq:TC}
| \Gamma(h_j^*) | \le \mbox{Tol} \qquad \mbox{and} \qquad |h_j^*-h_{j-1}^*| \le \mbox{TolR}|h_j^*|+ \mbox{TolA} \ ,
\end{equation}
with $ \mbox{Tol}=\mbox{TolR}=\mbox{TolA}=1\mbox{D}-06$.

%% file: CF2013TAB1.tex
\begin{table}[p]
\caption{Iterations for $\beta=-0.01$ with ${\ds \frac{d^2f^*}{d\eta^{*2}}(0)} = 1$.}
\vspace{.5cm}
\renewcommand\arraystretch{1.3}
	\centering
		\begin{tabular}{cr@{.}lr@{.}lcr@{.}l}
\hline 
{$j$} &
\multicolumn{2}{c}%
{$h_j^*$}
& \multicolumn{2}{c}%
{$\Gamma(h_j^*)$}
&
{${\ds \frac{|h_j^*-h_{j-1}^*|}{|h_j^*|}}$}
& \multicolumn{2}{c}%
{${\displaystyle \frac{d^2f}{d\eta^2}(0)}$} \\[1.2ex]
\hline
0 &  5 &            &     0 & 631459 & & 0 & 431723 \\
1 & 10 &            &      1 & 791425  & & 0 & 384034 \\
2 &  2 & 278111 & $-$0 & 182888 & 3.389602 & 0 & 454658 \\
3 &  2 & 993420 &     0 & 0465208 & 0.238960 & 0 & 454658 \\
4 &  2 & 848366 &     9 & 5$\mbox{D}-04$ & 0.050925 & 0 & 456418 \\
5 &  2 & 845340 & $-$5 & 0$\mbox{D}-06$ & 0.001064 & 0 & 456455 \\
6 &  2 & 845356 &     6 & 1$\mbox{D}-08$ & 5.6$\mbox{D}-06$ & 0 & 456455 \\
7 &  2 & 845355 &     7 & 3$\mbox{D}-10$ & 6.7$\mbox{D}-08$ & 0 & 456455 \\
\hline			
		\end{tabular}
	\label{tab:Itera1}
\end{table}

%% file: CF2013TAB2.tex
\begin{table}[p]
\caption{Iterations for $\beta=-0.01$ with ${\ds \frac{d^2f^*}{d\eta^{*2}}(0)} = -1$.}
\vspace{.5cm}
\renewcommand\arraystretch{1.3}
	\centering
		\begin{tabular}{cr@{.}lr@{.}lcr@{.}l}
\hline 
{$j$} &
\multicolumn{2}{c}%
{$h_j^*$}
& \multicolumn{2}{c}%
{$\Gamma(h_j^*)$}
&
{${\ds \frac{|h_j^*-h_{j-1}^*|}{|h_j^*|}}$}
& \multicolumn{2}{c}%
{${\displaystyle \frac{d^2f}{d\eta^2}(0)}$} \\[1.2ex]
\hline
0 &  75 &            &     0 & 731890 & & $-$0 & 059237 \\
1 & 150 &           &     5 & 263092 & & $-$0 & 092368 \\
2 &  62 & 885833 & $-$0 & 443040 & 1.385275 & $-$0 & 028870 \\
3 &  69 & 649620 &     0 & 181067 & 0.097112 & $-$0 & 046991 \\
4 & 67 & 687299 & $-$0 & 011297  & 0.028991 & $-$0 & 042016 \\
5 & 67 & 802542 & $-$2 & 1$\mbox{D}-04$ & 0.001700 & $-$0 & 042315 \\
6 & 67 & 804749 & 2 & 8$\mbox{D}-07$      & 3.3$\mbox{D}-05$ & $-$0 & 042321 \\
7 & 67 & 804746 & 7 & 9$\mbox{D}-10$    & 4.3$\mbox{D}-08$ & $-$0 & 042321 \\
\hline			
		\end{tabular}
	\label{tab:Itera2}
\end{table}

%% file: CF2013TAB4.tex
\begin{table}[!htb]
\caption{Comparison for the reverse flow skin-friction coefficients ${\displaystyle \frac{d^2f}{d\eta^{2}}(0)}$.
For all cases we used $h_0^* = 15$ and $h_1^* = 25$.
The iterations were, from top to bottom line: $8$, $7$, $9$, $7$, and $7$.}
\vspace{.5cm}
\renewcommand\arraystretch{1.3}
	\centering
		\begin{tabular}{r@{.}lcccc}
\hline 
\multicolumn{2}{c}%
{$\beta$} &
Stewartson \cite{Stewartson:1954:FSF} &
Asaithambi \cite{Asaithambi:1997:NMS} & Auteri et al. \cite{Auteri:2012:GLS}&
ITM \\[1ex]
\hline
$-0$ & $025$ & $-0.074$ & & $$ & $-0.074366$ \\
$-0$ & $05$  & $-0.108$ & & $$ & $-0.108271$ \\
$-0$ & $1$    & $-0.141$ & $-0.140546$ & $-0.140546$ & $-0.140546$ \\
$-0$ & $15$  & $-0.132$ & $-0.133421$ & $-0.133421$ & $-0.133421$ \\
$-0$ & $18$  & $-0.097$ & $-0.097692$ & $-0.097692$ & $-0.097692$ \\
\hline			
		\end{tabular}
	\label{tab:Compare}
\end{table}